\documentclass[11pt]{article}

\usepackage{latexsym,amssymb}

\makeatletter \@addtoreset{equation}{section}

\newtheorem{thm}{Theorem}[section]
\newtheorem{prop}[thm]{Proposition}
\newtheorem{lem}[thm]{Lemma}
\newtheorem{cor}[thm]{Corollary}

\newtheorem{example}[thm]{Example}

\def\enu#1{\newline\makebox[5mm][l]{\rm(#1)}}
\def\benu#1#2{\newline\makebox[#1 mm][l]{\rm(#2)}}

\def\bp{\noindent{\it Proof.}\ }
\def\bpp#1{\noindent{\it Proof of #1.}\ }
\def\ep{\nopagebreak\newline\mbox{\ }\hfill\rule{2mm}{2mm}}
\def\epp{\nopagebreak\mbox{\ }\hfill\rule{2mm}{2mm}}

\def\Ad{{\rm Ad}\,}

\def\det{{\rm det}}

\def\Irr{{\rm Irr}}
\def\Ker{{\rm Ker}}

\def\supp{{\rm supp}\,}
\def\Tr{{\rm Tr}}

\def\7#1{{\mathbb #1}}

\def\A{{\cal A}}

\def\C{{\cal C}}

\def\H{{\cal H}}

\def\iso{{\raisebox{1mm}{$\sim$}\hspace{-4.2mm}\to}}
\def\liso{{\raisebox{1mm}{$\sim$}\hspace{-4.6mm}\leftarrow}}

\def\2{{\frac{1}{2}}}
\def\ii{{\frac{i}{2}}}

\def\D{{\Delta}}
\def\sl2{{U_q({\mathfrak su}_2)}}

\def\eps{{\varepsilon}}

\begin{document}

\title{\bf Poisson boundary of the dual of $SU_q(n)$}

\author{Masaki Izumi$^1$, Sergey Neshveyev$^2$ and Lars Tuset$^3$}

\date{}

\footnotetext[1]{Supported by JSPS.}

\footnotetext[2]{Partially supported by the Norwegian Research
Council.}

\footnotetext[3]{Supported by the SUP-program of the Norwegian
Research Council.}

\maketitle

\begin{abstract}
We prove that for any non-trivial product-type action $\alpha$ of
$SU_q(n)$ ($0<q<1$) on an ITPFI factor $N$, the relative commutant
$(N^\alpha)'\cap N$ is isomorphic to the algebra
$L^\infty(SU_q(n)/\7T^{n-1})$ of bounded measurable functions on
the quantum flag manifold $SU_q(n)/\7T^{n-1}$. This is equivalent
to the computation of the Poisson boundary of the dual discrete
quantum group $\widehat{SU_q(n)}$. The proof relies on a
connection between the Poisson integral and the Berezin transform.
Our main technical result says that a sequence of Berezin
transforms defined by a random walk on the dominant weights of
$SU(n)$ converges to the identity on the quantum flag manifold.
This is a $q$-analogue of some known results on quantization of
coadjoint orbits of Lie groups.
\end{abstract}

\bigskip

\section*{Introduction}

The study of random walks on duals of compact groups, by which is
meant the study of convolution operators on group von Neumann
algebras of compact groups, was initiated by Biane in~\cite{B1}.
His results parallel the theory of random walks on discrete
abelian groups, and demonstrate some new interesting
non-commutative phenomena~\cite{B2,B3}. As was observed by the
first author~\cite{I}, the situation is even more interesting when
one considers compact quantum groups~\cite{W}. He was motivated by
the study of product-type actions of such groups on infinite
tensor products of factors of type~I (ITPFI). In the classical
case such actions are always minimal. For quantum groups this is
not so. In fact, the relative commutant of the fixed point algebra
is isomorphic to the algebra of bounded measurable functions on
the Poisson boundary of the dual discrete quantum group, that is,
to the algebra of bounded harmonic elements. The general theory
developed in~\cite{I} was illustrated by the computation of the
Poisson boundary of~$\widehat{SU_q(2)}$, which was shown to be
isomorphic to the quantum sphere~$S^2_q$. Later the second and the
third author~\cite{NT} computed the Martin boundary
of~$\widehat{SU_q(2)}$, that is, described all (unbounded)
harmonic functions, thus generalizing the results in~\cite{B3} and
establishing a connection between the results in~\cite{B3}
and~\cite{I}.

There are important differences between the classical boundary
theory (or even that for duals of compact groups as considered by
Biane) and the boundary theory for discrete quantum groups, which
we want to stress. In the classical setting, if one forgets about
the action of a group on its boundary, the computation of the
Poisson boundary reduces to the question whether it is trivial or
not. In the quantum case just the description of the
(noncommutative) algebra of functions on the Poisson boundary is a
highly nontrivial problem. On the other hand, such central
question of the classical theory as the description of minimal
harmonic functions, becomes of peripheral interest in the quantum
setting.

The purpose of the present paper is to compute the Poisson
boundary of~$\widehat{SU_q(n)}$. As we already said, this problem
can be formulated without any reference to random walks. Suppose
one is given a product-type action $\alpha$ of $SU_q(n)$ on~$N$.
By restriction we get an action of $SU_q(n)$ on~$(N^\alpha)'\cap
N$. If the action $\alpha$ is faithful, then $(N^\alpha)'\cap N$
is anti-isomorphic to $N'\cap(SU_q(n)\ltimes N)$. In this case the
dual action of $\widehat{SU_q(n)}$ on $N'\cap(SU_q(n)\ltimes N)$
induces an action on $(N^\alpha)'\cap N$. The problem is to
compute $(N^\alpha)'\cap N$ together with the two above actions of
$SU_q(n)$ and $\widehat{SU_q(n)}$. Our main result can be stated
in this setting as follows.

\medskip

\noindent {\bf Theorem A} {\it For any non-trivial product-type
action $\alpha$ of $SU_q(n)$ ($0<q<1$) on an ITPFI factor $N$, the
relative commutant $(N^\alpha)'\cap N$ is $SU_q(n)$-equivariantly
isomorphic to the algebra $L^\infty(SU_q(n)/\7T^{n-1})$ of bounded
measurable functions on the quantum flag manifold
$SU_q(n)/\7T^{n-1}$. When the action of $\widehat{SU_q(n)}$ on
$(N^\alpha)'\cap N$ is well-defined, the isomorphism is also
$\widehat{SU_q(n)}$-equivari\-ant.}

\medskip

As for the boundary interpretation of $(N^\alpha)'\cap N$, the
action of $\widehat{SU_q(n)}$ corresponds to the usual action by
translations of a group on its boundary, while the action of
$SU_q(n)$ can be thought of as coming from the symmetry group of
the measure. Then an equivalent form of Theorem~A is the
following.

\medskip

\noindent {\bf Theorem B} {\it Let $G$ be the $q$-deformation
($0<q<1$) of the compact group $SU(n)$, or of its quotient by a
normal subgroup. Let $T\subset G$ be the maximal torus. Then for
any $G$-invariant normal generating state on $l^\infty(\hat G)$,
the corresponding Poisson boundary of $\hat G$  is $G$- and $\hat
G$-equivariantly isomorphic to the quantum flag manifold $G/T$.}

\medskip

Note that even for $n=2$ this result is stronger than the ones
established in~\cite{I,NT}, since we make no additional
assumptions on the state.

There are several heuristic reasons why such a result should be
true. One is implied by the work of Biane~\cite{B3}. Although he
was able to complete his computations only for $SU(2)$
(see~\cite{C} for an extension of some of his results to~$SU(n)$),
his work strongly suggests that the Martin boundary of $\hat G$
for any compact semisimple Lie group~$G$, is the sphere in the
dual of the Lie algebra, and the action of $G$ is just the
coadjoint action. Thus one would expect that the Martin boundary
of $\widehat{SU_q(n)}$ is a certain quantization of the sphere. In
the classical case the Poisson boundary is the support of a
canonical measure on the Martin boundary. By~\cite{H} the action
of $SU_q(n)$ on the Poisson boundary of $\widehat{SU_q(n)}$ is
ergodic. This corresponds in the classical case to the fact that
the measure is supported on some orbit, and indeed, the typical
coadjoint orbit is the flag manifold.

Since, however, the Martin boundary seems difficult to compute, we
won't follow the approach suggested above. Our computation of the
Poisson boundary is based on a careful study of the completely
positive map $\Theta$ from $L^\infty(SU_q(n))$ into the algebra of
harmonic functions introduced in~\cite{I}. We call this map the
Poisson integral. We show that multiplicativity of $\Theta$
restricted to $L^\infty(SU_q(n)/\7T^{n-1})$ is equivalent to
convergence to the identity of a certain sequence of operators
on~$L^\infty(SU_q(n)/\7T^{n-1})$. These operators are analogues of
Berezin transforms. In the classical case it is known~\cite{D}
that Berezin transforms converge to the identity on the flag
manifold along any ray in the Weyl chamber. We prove that our
operators converge to the identity on the quantum flag manifold
along almost every path in the Weyl chamber. It is worth noticing
that though we benefited from these analogies, our proof bears no
relation to the classical proof. As a matter of fact, our
operators are not the most straightforward analogues of Berezin
transforms. In particular, in the classical limit they give
operators mapping everything to the scalars. Our proof of
convergence invokes yet another auxiliary operator, which in the
classical limit yields the identity operator, whereas in the
quantum case it is uniquely ergodic on~$C(SU_q(n)/\7T^{n-1})$.

Multiplicativity of $\Theta$ on $L^\infty(SU_q(n)/\7T^{n-1})$
implies injectivity. This can be interpreted as the existence of a
surjective map from the Poisson boundary onto the quantum flag
manifold. To show that this is an isomorphism we then use a
counting argument relying on the already mentioned ergodicity of
the action of $SU_q(n)$ on the Poisson boundary.

\bigskip

\section{Preliminaries} \label{s1}

We will use the same conventions for quantum groups as
in~\cite{NT}. We will, however, change the notation slightly to
hopefully make it more transparent. So a compact quantum group $G$
is defined by a unital C$^*$-algebra $C(G)$ with comultiplication
$\D\colon C(G)\to C(G)\otimes C(G)$, which is a unital
$*$-homomorphism such that $(\D\otimes\iota)\D=(\iota\otimes\D)\D$
and both $\D(C(G))(C(G)\otimes1)$ and $\D(C(G))(1\otimes C(G))$
are dense in $C(G)\otimes C(G)$. We will always work in the
reduced setting, that is, we assume that the Haar state $\varphi$
on $C(G)$ is faithful. In cases when it is more convenient to deal
with von Neumann algebras, we shall use the notation $L^\infty(G)$
for $\pi_\varphi(C(G))''$.

By a unitary representation of $G$ on a Hilbert space $H$ we mean
a unitary corepresentation of $(C(G),\D)$, that is, a unitary
$U\in M(C(G)\otimes B_0(H))$ such that
$(\D\otimes\iota)(U)=U_{13}U_{23}$. Here $B_0(H)$ denotes the
algebra of compact operators on $H$. Let $I=\Irr(G)$ be the set of
equivalence classes of irreducible unitary representations of $G$.
For each $s\in I$, we fix a representative $U^s\in C(G)\otimes
B(H_s)$. Then the algebra $c_0(\hat G)$ of functions on the dual
discrete quantum group $\hat G$ vanishing at infinity is defined
as the C$^*$-direct sum of $B(H_s)$, $s\in I$. The algebra
$l^\infty(\hat G)$ is defined as the W$^*$-direct sum of $B(H_s)$,
$s\in I$.

Let $\A(G)$ be the $*$-subalgebra of $C(G)$ generated by the
matrix coefficients of finite dimensional unitary representations
of $G$. Let also $\A(\hat G)\subset c_0(\hat G)$ be the algebraic
direct sum of $B(H_s)$, $s\in I$. There is a pairing between
$\A(G)$ and $\A(\hat G)$. If we fix matrix units
$\{m^s_{ij}\}_{i,j}$ in $B(H_s)$, and denote by
$\{u^s_{ij}\}_{i,j}$ the corresponding matrix coefficients of
$U^s$, then the pairing is given by
$(u^s_{ij},m^t_{kl})=\delta_{st}\delta_{ik}\delta_{jl}$. For a
unitary representation $U$ of $G$ on $H$, the corresponding
representation $\pi_U\colon l^\infty(\hat G)\to B(H)$ is given by
$$
\pi_U(\omega)=(\omega\otimes\iota)(U).
$$
In particular, if $W\in B(H_\varphi\otimes H_\varphi)$ is the
multiplicative unitary for $G$,
$$
W^*(\xi\otimes a\xi_\varphi)=\D(a)(\xi\otimes\xi_\varphi),
\ \ \xi\in H_\varphi,
\ \ a\in C(G),
$$
then we get a faithful representation $\pi_W\colon l^\infty(\hat
G)\to B(H_\varphi)$. Thus we can (and often will) think of
$L^\infty(G)$ and $l^\infty(\hat G)$ as subalgebras of
$B(H_\varphi)$. Note that the comultiplications $\D\colon
L^\infty(G)\to L^\infty(G)\otimes L^\infty(G)$ and $\hat\D\colon
l^\infty(\hat G)\to l^\infty(\hat G)\otimes l^\infty(\hat G)$ are
then given by
$$
\D(a)=W^*(1\otimes a)W,\ \ \hat\D(x)=W(x\otimes1)W^*.
$$

Given a unitary representation $U$ of $G$ on $H$, we define a left
and a right action of $G$ on $B(H)$ by
$$
\alpha_{U,l}\colon B(H)\to L^\infty(G)\otimes B(H),
\ \ \alpha_{U,l}(x)=U^*(1\otimes x)U,
$$
$$
\alpha_{U,r}\colon B(H)\to B(H)\otimes L^\infty(G),\ \
\alpha_{U,r}(x)=U_{21}(x\otimes1)U^*_{21}.
$$
As $l^\infty(\hat G)=\oplus_s B(H_s)$, we get in particular a left
and a right action of $G$ on $l^\infty(\hat G)$, which we denote
by $\Phi_l$ and $\Phi_r$, respectively. Thinking of $l^\infty(\hat
G)$ as the group von Neumann algebra of $G$, they are analogues of
the adjoint action. If the representation $U$ is finite
dimensional, then the actions of $G$ on $B(H)$ have canonical
invariant states
$$
\phi_U=\frac{1}{d_U}\Tr(\cdot\pi_U(\rho^{-1})), \ \
(\iota\otimes\phi_U)\alpha_{U,l}=\phi_U(\cdot)1,
$$
$$
\omega_U=\frac{1}{d_U}\Tr(\cdot\pi_U(\rho)), \ \
(\omega_U\otimes\iota)\alpha_{U,r}=\omega_U(\cdot)1,
$$
where
$$
\rho=f_1
$$
is the Woronowicz character and
$d_U=\Tr(\pi_U(\rho^{-1}))=\Tr(\pi_U(\rho))$ is the quantum
dimension of~$U$. We will write $\phi_s$ and $\omega_s$ instead of
$\phi_{U^s}$ and $\omega_{U^s}$, respectively. In general, all
normal left (resp. right) $G$-invariant states on $B(H)$ are given
by $\Tr(\cdot a\pi_U(\rho^{-1}))$ (resp. by  $\Tr(\cdot
a\pi_U(\rho))$),  where $a\in\pi_U(l^\infty(\hat G))'$ is a
positive element such that $\Tr(a\pi_U(\rho^{-1}))=1$ (resp.
$\Tr(a\pi_U(\rho))=1$).

Denote by $\C_l(\hat G)$ (resp. by $\C_r(\hat G)$) the space of
normal left (resp. right) $G$-invariant functionals on
$l^\infty(\hat G)$. Then $\C_l(\hat G)$ (resp. $\C_r(\hat G)$) is
the closed linear span of $\phi_s$ (resp. $\omega_s$), $s\in I$.
The space $\C_l(\hat G)$ (as well as $\C_r(\hat G)$) is an algebra
with product $\phi_1\phi_2=(\phi_1\otimes\phi_2)\hat\D$.
Alternatively, one can define a fusion algebra structure on
$R(G)=\oplus_{s\in I}\7Z$ in the sense of~\cite{HI}. Then the
quantum dimension function on $R(G)$ defines a convolution algebra
$l^1(I)$, which is isomorphic to $\C_l(\hat G)$ and $\C_r(\hat
G)$. The algebra $\C_l(\hat G)$ has an anti-linear involution
$\phi\mapsto\check\phi$ such that $\check\phi_U=\phi_{\bar U}$,
where $\overline{U}$ is the conjugate unitary representation. Let
also $s\mapsto\bar s$ be the involution on $I$, so $U^{\bar
s}\cong\overline{U^s}$.

\medskip

Given a normal state $\phi\in l^\infty(\hat G)_*$, define the
convolution operator $P_\phi=(\phi\otimes\iota)\hat\D$ on
$l^\infty(\hat G)$, and set
$$
H^\infty(\hat G,\phi)=\{x\in l^\infty(\hat G)\,|\, P_\phi(x)=x\}.
$$
Then $H^\infty(\hat G,\phi)$ is (the algebra of bounded measurable
functions on) the Poisson boundary of $\hat G$ with respect to
$\phi$. According to~\cite{I} the algebra structure on
$H^\infty(\hat G,\phi)$ is given by
$$
x\cdot y=s^*-\lim_{n\to\infty}\frac{1}{n}\sum^{n-1}_{k=0}
P^k_\phi(xy).
$$
We will only be interested in the case when $\phi\in\C_l(\hat G)$.
Denote by $\supp\phi\subset I$ the set of $s\in I$ such that
$\phi(I_s)>0$, where $I_s$ is the unit in $B(H_s)\subset
l^\infty(\hat G)$. A state $\phi$ is called generating if
$\cup_{n\in\7N}\supp \phi^n=I$. If $\phi$ is not generating, but
$\supp\phi$ is symmetric, then the norm closure of the linear span
of the matrix coefficients $u^s_{ij}$ for $s\in\cup_{n\in\7N}\supp
\phi^n$ is the algebra $C(H)$ for a compact quantum group $H$.
Thus $H$ is a quotient of $G$, and by the orthogonality relations
we have $\Irr(H)=\cup_{n\in\7N}\supp \phi^n$. The symmetry
assumption for $\supp\phi$ is needed to ensure that $C(H)$ is a
$*$-algebra. In the case of $q$-deformations of compact connected
semisimple Lie groups~\cite{KoS} this assumption is redundant.
Indeed, it is well-known that any compact connected semisimple Lie
group $G$ has the property that if $U$ is a faithful unitary
representation of $G$, then any irreducible representation of $G$
appears as a subrepresentation of the tensor power $U^{\times n}$
of $U$ for some $n\in\7N$. It follows that $\cup_{n\in\7N}\supp
\phi^n$ is always symmetric, more precisely, $\cup_{n\in\7N}\supp
\phi^n=\Irr(H)$ with $H=G/(\cap_{s\in\supp\phi}\Ker\,U^s)$. As the
fusion algebra is independent of the deformation parameter, we
conclude that for the $q$-deformation $G$ of any compact connected
semisimple Lie group, the set $\cup_{n\in\7N}\supp \phi^n$ is
symmetric, so it corresponds to a certain quotient of~$G$.

\medskip

Recall now the connection between Poisson boundaries and
product-type actions~\cite{I}. Let $U$ be a unitary representation
of $G$ on $H$, and $\tilde\phi\in B(H)_*$ a normal faithful
$\alpha_{U,l}$-invariant state. Set
$$
(N,\nu)=\mathop{\otimes}^{-1}_{-\infty}(B(H),\tilde\phi),
$$
and $N_n=\ldots\otimes1\otimes B(H)^{\otimes n}\subset N$. The
actions $\alpha_{U^{\times n},l}$ of $G$ on $B(H^{\otimes n})\cong
N_n$ define a left action $\alpha$ of $G$ on $N$, where $U^{\times
n}=U_{12}\ldots U_{1,n+1}$. To describe the relative commutant
$(N^\alpha)'\cap N$, consider $\phi=\tilde\phi\pi_U\in\C_l(\hat
G)$. Let $E_n\colon N\to N_n$ be the $\nu$-preserving conditional
expectation, and $j_n\colon l^\infty(\hat G)\to N$ the
homomorphism defined by
$j_n(x)=\ldots\otimes1\otimes\pi_{U^{\times n}}(x)$. Then
$$
E_nj_{n+1}=j_nP_\phi.
$$
The kernel of $j_n$ is $\oplus_{s\notin\supp\phi^n}B(H_s)$. So if
$F_n\colon l^\infty(\hat G)\to\oplus_{s\in\supp\phi^n}B(H_s)$ is
the canonical projection, we have $j_nF_n=j_n$, and $j_n$ is
faithful on $F_n(l^\infty(\hat G))$. As the image of
$(N^\alpha)'\cap N$ under $E_n$ is contained in the relative
commutant
$$
(B(H^{\otimes n})^{\alpha_{U^{\times n},l}})'\cap B(H^{\otimes n}),
$$
which coincides with the image of
$\pi_{U^{\times n}}$, we conclude that the elements of
$(N^\alpha)'\cap N$ are in one-to-one correspondence with bounded
sequences $\{x_n\}^\infty_{n=1}$ such that $x_n\in
F_n(l^\infty(\hat G))$ and
$$
x_n=F_nP_\phi(x_{n+1}).
$$
Such sequences are called $P_\phi$-harmonic. It turns out that if
$\phi$ is generating, then any bounded harmonic sequence
$\{x_n\}_n$ defines an element $x\in H^\infty(\hat G,\phi)$ such
that $x_n=F_n(x)$, and vice versa. Thus we get an isomorphism
$$
j_\infty\colon H^\infty(\hat G,\phi)\iso(N^\alpha)'\cap N,\ \
j_\infty(x)=s^*-\lim_{n\to\infty}j_n(x).
$$

By restricting the left action $\Phi_l$ of $G$ and the right
action $\hat\D$ of $\hat G$ on $l^\infty(\hat G)$ to
$H^\infty(\hat G,\phi)$, we get actions of $G$ and $\hat G$ on the
Poisson boundary. When $\phi$ is generating and we identify
$H^\infty(\hat G,\phi)$ with $(N^\alpha)'\cap N$, the action of
$G$ is just the restriction of $\alpha$ to the relative commutant,
since the homomorphisms $j_n$ are obviously $G$-equivariant. As
was conjectured in~\cite{I}, the action of $\hat G$ is related to
the dual action of $\hat G$ on $N'\cap(G\ltimes N)$. In the rest
of this section we want to clarify this point. This will not be
used in the subsequent sections.

\medskip

Consider first a more general situation. Suppose we are given a
left action of $G$ on a von Neumann algebra $N$,
$$
\alpha\colon N\to L^\infty(G)\otimes N.
$$
The crossed product $G\ltimes N$ is by definition the
W$^*$-subalgebra of $B(H_\varphi)\otimes N$ generated by
$l^\infty(\hat G)\otimes 1$ and $\alpha(N)$. Let $\nu$ be a normal
faithful $G$-invariant state on $N$. Then the map
$$
v\colon H_\nu\to H_\varphi\otimes H_\nu,\ \
a\xi_\nu\mapsto\alpha(a)(\xi_\varphi\otimes\xi_\nu)\ \
\hbox{for}\ a\in N,
$$
is an isometry. As $\alpha(N)(l^\infty(\hat G)\otimes 1)$ is dense
in $G\ltimes N$ and $x\xi_\varphi=\hat\eps(x)\xi_\varphi$ for
$x\in l^\infty(\hat G)$, where $\hat\eps$ is the counit on
$l^\infty(\hat G)$, we see that the image of $v$ is a $(G\ltimes
N)$-invariant subspace of $H_\varphi\otimes H_\nu$. Let $V$ be the
canonical representation of $G$ on $H_\nu$ implementing the
action,
$$
V^*(\xi\otimes a\xi_\nu)=\alpha(a)(\xi\otimes\xi_\nu),
\ \xi\in H_\varphi,\ a\in N.
$$
Since
\begin{eqnarray*}
(1\otimes v)V^*(\xi\otimes a \xi_\nu)
&=&(1\otimes v)\alpha(a)(\xi\otimes\xi_\nu)
=(\iota\otimes\alpha)\alpha(a)
(\xi\otimes\xi_\varphi\otimes\xi_\nu)\\
&=&(\D\otimes\iota)\alpha(a)(\xi\otimes\xi_\varphi\otimes\xi_\nu)
=(W^*\otimes1)(1\otimes\alpha(a))
(\xi\otimes\xi_\varphi\otimes\xi_\nu)\\
&=&(W^*\otimes1)(1\otimes v)(\xi\otimes a\xi_\nu),
\end{eqnarray*}
$V^*=(S\otimes\iota)(V)$ and $W^*=(S\otimes\iota)(W)$, where $S$
is the coinverse for $G$, we get
$$
v\pi_V(x)=(x\otimes1)v\ \ \hbox{for}\ x\in l^\infty(\hat G).
$$
We also have $v^*\alpha(a)v=a$ for $a\in N$. Thus
\begin{equation} \label{e1.10}
v^*(x\otimes1)\alpha(a)v=\pi_V(x)a\ \ \hbox{for}\
x\in l^\infty(\hat G),\ a\in N.
\end{equation}
Noting that $(J_\varphi\otimes J_\nu)v=vJ_\nu$, and recalling that
$J_\varphi x^*J_\varphi=\hat R(x)$ for $x\in l^\infty(\hat G)$,
where $\hat R$ is the unitary antipode on $l^\infty(\hat G)$, we
see that $J_\nu\pi_V(l^\infty(\hat G))J_\nu=\pi_V(l^\infty(\hat
G))$. Thus the map $x\mapsto J_\nu x^*J_\nu$ defines an
anti-isomorphism of $N'\cap(N\vee\pi_V(l^\infty(\hat G)))$ onto
$$
N\cap(N'\vee\pi_V(l^\infty(\hat G)))
=N\cap(N\cap\pi_V(l^\infty(\hat G))')'=N\cap(N^\alpha)'.
$$
It follows that the map
$$
\theta\colon N'\cap(G\ltimes N)\to (N^\alpha)'\cap N,
\ \ x\mapsto J_\nu v^*x^*vJ_\nu,
$$
is a normal surjective $*$-anti-homomorphism.

\smallskip

The dual action $\hat\alpha\colon G\ltimes N\to l^\infty(\hat
G)\otimes (G\ltimes N)$ is defined by
$$
\hat\alpha(x)=\tilde W(1\otimes x)\tilde W^*, \ \
x\in G\ltimes N,
$$
where
$\tilde W=((J_\varphi\otimes J_\varphi)W_{21}
(J_\varphi\otimes J_\varphi))\otimes1$,
so that
$$
\hat\alpha((x\otimes1)\alpha(a))
=(\hat\D(x)\otimes1)(1\otimes\alpha(a))
\ \ \hbox{for}\ a\in N,\ x\in l^\infty(\hat G).
$$
The left action $\hat\alpha$ of $\hat G$ on $G\ltimes N$ induces a
right action $\hat\alpha^{op}$ of $\hat G$ on $(G\ltimes N)^{op}$,
$$
\hat\alpha^{op}(x)=(\iota\otimes\hat R)(\alpha(x)_{21}).
$$

\smallskip

In the simplest case when $N=B(H)$ all these constructions can be
made more explicit.

\begin{lem} \label{1.1}
Let $U$ be a unitary representation of $G$ on $H$, $N=B(H)$,
$\alpha=\alpha_{U,l}$, $\nu$ a normal faithful $G$-invariant state
on $B(H)$, $\phi=\nu\pi_U$. Then the map
$$
\theta_0\colon l^\infty(\hat G)\otimes N\to G\ltimes N, \ \
\theta_0(x)=U^*xU,
$$
is an isomorphism with the following properties: \enu{i}
$(\iota\otimes\theta_0)(\hat\D\otimes\iota)(x)
=\hat\alpha\theta_0(x)$ for $x\in l^\infty(\hat G)\otimes N$;
\enu{ii} for the conditional expectation
$E_0=\iota\otimes\nu\colon G\ltimes N \to l^\infty(\hat G)$, we
have $E_0\theta_0(x\otimes1)=\hat RP_\phi\hat R(x)$ for $x\in
l^\infty(\hat G)$; \enu{iii} $\theta\theta_0(x\otimes1)=\pi_U\hat
R(x)$ for $x\in l^\infty(\hat G)$.
\end{lem}

\bp Since $(\iota\otimes\pi_U)\hat\D(l^\infty(\hat G))$ and
$1\otimes N$ generate $l^\infty(\hat G)\otimes N$, the fact that
$\theta_0$ is an isomorphism follows from the identities
$$
U^*(1\otimes a)U=\alpha(a)\ \ \hbox{for}\ a\in N,\ \
U^*(\iota\otimes\pi_U)\hat\D(x)U=x\otimes1\ \ \hbox{for}\ x\in
l^\infty(\hat G),
$$
where we have used that $U=(\iota\otimes\pi_U)(W)$.

\smallskip

The equality in (i) obviously holds for $x=1\otimes a$, $a\in N$.
Thus it is enough to consider $x=(\iota\otimes\pi_U)\hat\D(y)$.
Then $\theta_0(x)=y\otimes1$, so
$\hat\alpha\theta_0(x)=\hat\D(y)\otimes1$. On the other hand,
\begin{eqnarray*}
(\iota\otimes\theta_0)(\hat\D\otimes\iota)(x)
&=&(\iota\otimes\theta_0)(\iota\otimes\iota\otimes\pi_U)
(\hat\D\otimes\iota)\hat\D(y)\\
&=&(\iota\otimes\theta_0)(\iota\otimes\iota\otimes\pi_U)
(\iota\otimes\hat\D)\hat\D(y)
=\hat\D(y)\otimes1.
\end{eqnarray*}
Thus (i) is proved.

\smallskip

To prove (ii), note that for any $y,z\in l^\infty(\hat G)$ we have
$$
E_0\theta_0(\iota\otimes\pi_U)(\hat\D(y)(1\otimes z))
=E_0((y\otimes1)\alpha\pi_U(z))=\nu\pi_U(z)y=\phi(z)y.
$$
In other words, if we define a map $r\colon \A(\hat
G)\otimes\A(\hat G) \to \A(\hat G)\otimes \A(\hat G)$ by
$r(y\otimes z)=\hat\D(y)(1\otimes z)$, then
$$
E_0\theta_0(\iota\otimes\pi_U)r=\iota\otimes\phi.
$$
It is well-known (see e.g. \cite[Proposition 1.3]{NT}) that the
map $r$ is bijective with inverse $s$ given by $s(y\otimes
z)=(\iota\otimes\hat S)((1\otimes \hat S^{-1}(z))\hat\D(y))$.
Hence for any $y,z\in\A(\hat G)$ we get
$$
E_0\theta_0(y\otimes\pi_U(z))=(\iota\otimes\phi)s(y\otimes z)
=(\iota\otimes\phi\hat S)((1\otimes \hat S^{-1}(z))\hat\D(y)).
$$
As the state $\phi$ is $G$-invariant and $\hat S=\rho^\2\hat
R(\cdot)\rho^{-\2}$, we have $\phi\hat S=\phi\hat R$. Choosing a
net $\{z_i\}_i\subset\A(\hat G)$ of central projections converging
strongly to the unit, we thus get
$$
E_0\theta_0(y\otimes1)=(\iota\otimes\phi\hat R)\hat\D(y)
=\hat R(\iota\otimes\phi)(\hat R\otimes\hat R)\hat\D(y)
=\hat R(\phi\otimes\iota)\hat\D\hat R(y)=\hat RP_\phi\hat R(y),
$$
where we have used that $(\hat R\otimes\hat
R)\hat\D=\hat\D^{op}\hat R$. This proves (ii).

\smallskip

To prove (iii) we identify $H_\nu$ with the space $HS(H)$ of
Hilbert-Schmidt operators on $H$, so that $\xi_\nu=T^\2_0$ if
$T_0\in B(H)$ is the operator defining $\nu$, $\nu=\Tr(\cdot
T_0)$. If we further identify $HS(H)$ with $\overline{H}\otimes
H$, we have $J_\nu(\bar\xi\otimes\zeta)=\bar\zeta\otimes\xi$, and
$a(\bar\xi\otimes\zeta)=\bar\xi\otimes a\zeta$ for $a\in N=B(H)$.
Thus we have to prove that
$\theta\theta_0(x\otimes1)=1\otimes\pi_U\hat R(x)\in
B(\overline{H}\otimes H)$ for $x\in l^\infty(\hat G)$.

To compute $\theta$ we will check first that $V=\overline{U}\times
U$. The unitary $V\in M(C(G)\otimes B_0(H_\nu))$ is characterized
by the properties that it implements the action, $V^*(1\otimes
a)V=\alpha(a)$ for $a\in N$, and that
$V(\xi\otimes\xi_\nu)=\xi\otimes\xi_\nu$ for $\xi\in H_\varphi$.
The first property is obviously satisfied by $\overline{U}\times
U$. To check the second one, recall that by definition
$\overline{U}=(R\otimes j)(U)\in M(C(G)\otimes
B_0(\overline{H}))$, where $R$ is the unitary antipode on $C(G)$
and $j(x)\bar\xi=\overline{x^*\xi}$. If we consider $C(G)\otimes
HS(H)$ as a left module over $M(C(G)\otimes
B_0(\overline{H}\otimes H))$, it is then easy to check that
$$
(\overline{U}\times U)(1\otimes T)
=(R\otimes\iota)((R\otimes\iota)(U)(1\otimes T)U)\ \
\hbox{for}\ T\in HS(H).
$$
We have to prove that $(\overline{U}\times U)(1\otimes
T^\2_0)=1\otimes T^\2_0$. This follows from the identities
$$
(R\otimes\iota)(U)
=(1\otimes\pi_U(\rho^{-\2}))(S\otimes\iota)(U)
(1\otimes\pi_U(\rho^\2))
=(1\otimes\pi_U(\rho^{-\2}))U^*
(1\otimes\pi_U(\rho^\2))
$$
and
$$
(1\otimes T^\2_0)U=(1\otimes\pi_U(\rho^{-\2}))U
(1\otimes\pi_U(\rho^\2))
(1\otimes T^\2_0),
$$
where we have used  that $\pi_U(\rho)T_0\in \pi_U(l^\infty(\hat
G))'$. Thus $V=\overline{U}\times U$.

By virtue of (\ref{e1.10}), for any $y,z\in l^\infty(\hat G)$ we
then have
\begin{eqnarray*}
v^*\theta_0((\iota\otimes\pi_U)(\hat\D(y)(1\otimes z)))v
&=&v^*(y\otimes1)\alpha\pi_U(z)v
=\pi_V(y)(1\otimes\pi_U(z))\\
&=&(\pi_{\bar U}\otimes\pi_U)\hat\D(y)(1\otimes\pi_U(z))
=(\pi_{\bar U}\otimes\pi_U)(\hat\D(y)(1\otimes z)).
\end{eqnarray*}
Since $\hat\D(l^\infty(\hat G))(1\otimes l^\infty(\hat G))$ is
dense in $l^\infty(\hat G)\otimes l^\infty(\hat G)$, for any $x\in
l^\infty(\hat G)$ we thus get
$$
v^*\theta_0(x\otimes1)v=(\pi_{\bar U}\otimes\pi_U)(x\otimes1)
=\pi_{\bar U}(x)\otimes1,
$$
whence
$$
\theta\theta_0(x\otimes 1)=J_\nu v^*\theta_0(x^*\otimes1)vJ_\nu
=J_\nu(\pi_{\bar U}(x^*)\otimes1)J_\nu
=1\otimes\overline{\pi_{\bar U}(x^*)},
$$
where $\overline{\pi_{\bar U}(x^*)}\xi=\overline{\pi_{\bar
U}(x^*)\bar\xi}$. By definition of the conjugate representation we
have
$$
\pi_{\bar U}(x^*)\bar\xi=\overline{\pi_U\hat R(x)\xi}.
$$
It follows that $\theta\theta_0(x\otimes 1)=1\otimes\pi_U\hat
R(x)$. \ep

Return now to the study of $N'\cap(G\ltimes N)$ for a product-type
action defined by a representation $U$ of $G$ on $H$ and a
$G$-invariant state $\tilde\phi\in B(H)_*$. The relative commutant
can be described in a way similar to that for $(N^\alpha)'\cap N$,
see \cite{V}. The conditional expectation $E_n\colon N\to N_n$
extends to the conditional expectation $\iota\otimes E_n\colon
G\ltimes N\to G\ltimes N_n$. Let $y\in N'\cap(G\ltimes N)$. Then
$(\iota\otimes E_n)(y)\in N_n'\cap(G\ltimes N_n)$, so by
Lemma~\ref{1.1}
$$
(\iota\otimes E_n)(y)=(U^{\times n})^*(y_n\otimes1)U^{\times n}
$$
for a unique element $y_n\in l^\infty(\hat G)$. As
$E_nE_{n+1}=E_n$, we have
\begin{eqnarray*}
(U^{\times n})^*(y_n\otimes1)U^{\times n}
&=&(\iota\otimes E_n)\Big((U^{\times (n+1)})^*(y_{n+1}\otimes1)
U^{\times(n+1)}\Big)\\
&=&(U^{\times n})^*((\iota\otimes\tilde\phi)
(U^*(y_{n+1}\otimes1)U)\otimes1)U^{\times n},
\end{eqnarray*}
whence by Lemma~\ref{1.1}(ii)
$$
y_n=(\iota\otimes\tilde\phi)(U^*(y_{n+1}\otimes1)U)
=\hat RP_\phi\hat R(y_{n+1}),
$$
where $\phi=\tilde\phi\pi_U$. Taking into account parts (i) and
(iii) of Lemma~\ref{1.1} as well as the description of
$(N^\alpha)'\cap N$ in terms of harmonic sequences, we get the
following result.

\begin{prop}
For the product-type action $\alpha$ of a compact quantum group
$G$ defined by a representation $U$ of $G$ on $H$ and a
$G$-invariant normal faithful state $\tilde\phi$ on $B(H)$, there
exists a $\hat G$-equivariant complete order isometry between
$N'\cap(G\ltimes N)$ and the space of bounded sequences
$\{y_n\}^\infty_{n=1}\subset l^\infty(\hat G)$ such that $\hat
R(y_n)=P_\phi\hat R(y_{n+1})$, where $\phi=\tilde\phi\pi_U$.

The surjective anti-homomorphism $\theta\colon N'\cap(G\ltimes
N)\to(N^\alpha)'\cap N $ maps the element defined by a sequence
$\{y_n\}^\infty_{n=1}$ to the element defined by the
$P_\phi$-harmonic sequence $\{F_n\hat R(y_n)\}^\infty_{n=1}$.\epp
\end{prop}

Note that even if $\phi$ is generating, the homomorphism
$(N'\cap(G\ltimes N))^{op}\to(N^\alpha)'\cap N$ needs not be
injective nor $\hat G$-equivariant (where we consider
$(N'\cap(G\ltimes N))^{op}$ with the right action
$\hat\alpha^{op}$ of $\hat G$, and identify $(N^\alpha)'\cap N$
with $H^\infty(\hat G,\phi)$ to get a right action of $\hat G$ on
$(N^\alpha)'\cap N$).

\begin{example}\rm
Let $G=SU(2)$ and $U$ be the fundamental representation of $G$.
Identify $\Irr(G)$ with $\2\7Z_+$ as usual. Set
$p=\sum^\infty_{n=1}I_{n-\2}$. Then
$$
\hat\D(p)=p\otimes(1-p)+(1-p)\otimes p.
$$
In particular, $P_\phi(p)=1-p$. Consider the sequence
$\{y_n\}^\infty_{n=1}$ defined by $y_{2k+1}=1-p$, $y_{2k}=p$, and
let $y$ be the corresponding element of $N'\cap(SU(2)\ltimes N)$.
Then $y$ is in the kernel of $\theta$. Moreover, as
$$
\hat\alpha(y)=p\otimes(1-y)+(1-p)\otimes y,
$$
the $\hat G$-equivariance also does not hold.
\end{example}

It is clear, however, what goes wrong in the previous example.
Though the representation $U$ was faithful, $\overline{U}\times U$
was not, so we had in fact an action of $SO(3)$. More generally,
assume $G$ is the $q$-deformation of a compact connected
semisimple Lie group. Let $H$ be the quotient of $G$ defined by
$\overline{U}\times U$, that is, $C(H)$ is the C$^*$-subalgebra of
$C(G)$ generated by the matrix coefficients of $\overline{U}\times
U$. Then $\Irr(H)=\cup_{n\in\7N}\supp(\check\phi\phi)^n$. Choose
also $k\in\7N$ such that $0\in\supp\phi^k$, that is, $U^{\times
k}$ contains the trivial representation. Since $0\in\supp\phi^k$,
we have $\supp\phi^k\subset\supp(\check\phi^k\phi^k)
=\supp(\check\phi\phi)^k\subset\Irr(H)$. We have thus shown that
by replacing $G$ by its quotient and $U$ by some power $U^{\times
k}$, in the study of product-type actions we can always assume
that $\overline{U}\times U$ is faithful, that is, $C(G)$ is
generated by the matrix coefficients of $\overline{U}\times U$.

Assume now that $\overline{U}\times U$ is faithful. Then
$\Irr(G)=\cup_{n\in\7N}\supp\phi^{kn}$ for any $k\in\7N$. Indeed,
as we discussed earlier, the set $\cup_{n\in\7N}\supp\phi^{kn}$ is
always symmetric. Hence it contains $\supp(\check\phi\phi)^{kn}$
for any $n\in\7N$. Since $0\in\supp(\check\phi\phi)$, we have
$\supp(\check\phi\phi)^m \subset\supp(\check\phi\phi)^{m+1}$ for
any $m$. Thus
$$
\Irr(G)=\cup_{n\in\7N}\supp(\check\phi\phi)^n
=\cup_{n\in\7N}\supp(\check\phi\phi)^{kn}
\subset\cup_{n\in\7N}\supp\phi^{kn}.
$$
In particular, we can find $k$ and $n$ such that $0\in\supp\phi^k$
and $\supp\phi\cap\supp\phi^{kn}\ne\emptyset$. As
$0\in\supp\phi^{kn}$, we have
$$
\supp\phi\cap\supp\phi^{kn}\subset\supp\phi^{kn}\cap\supp\phi^{kn+1}.
$$
By the $0$-$2$ law (see e.g \cite[Proposition 2.12]{NT}) we can
conclude that $\|P^m_\phi-P^{m+1}_\phi\|\to0$ as $m\to\infty$. It
follows that if $\{x_n\}^\infty_{n=1}$ is a bounded sequence in
$l^\infty(\hat G)$ such that $x_n=P_\phi(x_{n+1})$, then the
sequence is constant, $x_n=x_{n+1}$. Note also that if $y\in
N'\cap(G\ltimes N)$ is the element defined by the sequence $\{\hat
R(x_n)\}^\infty_{n=1}$, then by Lemma~\ref{1.1}(ii), for the
conditional expectation $E_0=\iota\otimes\nu\colon G\ltimes N\to
l^\infty(\hat G)$ we have $E_0(y)=\hat RP_\phi(x_1)$. We have thus
proved the following result.

\begin{cor}
Let $G$ be the $q$-deformation of a compact connected semisimple
Lie group, $\alpha$ the product-type action of $G$ on $(N,\nu)$
defined by a representation $U$ of $G$ on $H$ and a $G$-invariant
normal faithful state $\tilde\phi$ on $B(H)$. Assume that the
representation $\overline{U}\times U$ is faithful. Then we have
isomorphisms
$$
(N'\cap(G\ltimes N))^{op}\mathop{\iso}_\theta(N^\alpha)'\cap N
\mathop{\liso}_{j_\infty}H^\infty(\hat G,\phi),
$$
where $\phi=\tilde\phi\pi_U$. Moreover, the homomorphism
$j^{-1}_\infty\theta$ is $\hat G$-equivariant and coincides with
the map $\hat R\otimes\nu$.\epp
\end{cor}

For a more general group $G$ and a generating state
$\phi\in\C_l(\hat G)$, the homomorphism $j^{-1}_\infty\theta$ is a
$\hat G$-equivariant isomorphism if any bounded sequence
$\{x_n\}^\infty_{n=1}\subset l^\infty(\hat G)$ such that
$P_\phi(x_{n+1})=x_n$ is constant, that is, $x_{n+1}=x_n$. For
this it is enough to require that
$\supp\phi^m\cap\supp\phi^{m+1}\ne\emptyset$ for some $m$.

\bigskip

\section{Poisson integral and Berezin transform} \label{s2}

Let $G$ be a compact quantum group and $\phi\in\C_l(\hat G)$ be a
normal left $G$-invariant state on $l^\infty(\hat G)$.
Consider the right action $\hat\Phi$ of
$\hat G$ on $L^\infty(G)$,
$$
\hat\Phi\colon L^\infty(G)\to L^\infty(G)\otimes l^\infty(\hat G),
\ \ \hat\Phi(a)=W(a\otimes1)W^*.
$$
It was proved in \cite{I} that
$$
\Theta=(\varphi\otimes\iota)\hat\Phi
$$
maps $L^\infty(G)$ into $H^\infty(\hat G,\phi)$. In fact, this is
the only normal unital $G$- and $\hat G$-equivariant map from
$L^\infty(G)$ into $l^\infty(\hat G)$ (we consider $L^\infty(G)$
and $l^\infty(\hat G)$ with the left actions of $G$ given by $\D$
and $\Phi_l$, respectively, and with the right actions of $\hat G$
given by $\hat\Phi$ and $\hat\D$, respectively). Indeed, assume
$T$ is such a map. Since the counit $\hat\eps$ on $l^\infty(\hat
G)$ is $G$-invariant, $\hat\eps T$ is a $G$-invariant normal
linear functional on $L^\infty(G)$ such that $\hat\eps T(1)=1$.
Hence $\hat\eps T=\varphi$. Then by $\hat G$-equivariance of $T$,
we get
$$
T=(\hat\eps\otimes\iota)\hat\D T
=(\hat\eps\otimes\iota)(T\otimes\iota)\hat\Phi
=(\varphi\otimes\iota)\hat\Phi=\Theta.
$$
Thus if the Poisson boundary is $G$- and $\hat G$-equivariantly
isomorphic to a homogeneous space $G/H$ of $G$, the map $\Theta$
must be an isomorphism of $L^\infty(G/H)$ onto $H^\infty(\hat
G,\phi)$. We call the map $\Theta$ the Poisson integral.

To show multiplicativity of $\Theta$ we will use the following
criterion.

\begin{lem} \label{2.1}
Let $N_i$ be a von Neumann algebra, $\nu_i$ a normal faithful
state on $N_i$, $i=1,2$, $\theta\colon N_1\to N_2$ a normal ucp
map such that $\nu_2\theta=\nu_1$ and
$\sigma^{\nu_2}_t\theta=\theta\sigma^{\nu_1}_t$. Then there exists
a normal ucp map $\theta^*\colon N_2\to N_1$ such that
$$
\nu_2(\theta(x_1)x_2)=\nu_1(x_1\theta^*(x_2))
\ \ \hbox{for}\ x_1\in N_1,\ x_2\in N_2.
$$
For any $x\in N_1$, the following conditions are then equivalent:
\benu{6}{i} $x$ is in the multiplicative domain of $\theta$;
\benu{6}{ii} $||\theta(x)||_2=||x||_2$; \benu{6}{iii}
$||\theta(x^*)||_2=||x^*||_2$; \benu{6}{iv} $\theta^*\theta(x)=x$.
\end{lem}

\bp Since $\sigma^{\nu_2}_t\theta=\theta\sigma^{\nu_1}_t$, the map
$\theta^*$ can equivalently be defined by the condition
$$
(\theta(x_1)J_{\nu_2}x^*_2\xi_{\nu_2},\xi_{\nu_2})
=\nu_2(\theta(x_1)\sigma^{\nu_2}_{-\ii}(x_2))
=\nu_1(x_1\sigma^{\nu_1}_{-\ii}\theta^*(x_2))
=(x_1J_{\nu_1}\theta^*(x^*_2)\xi_{\nu_1},\xi_{\nu_1}).
$$
Hence $\theta^*$ exists \cite{AC}. A more general result than the
equivalence of the conditions (i)-(iv) can be found in~\cite{P}.
We will give a proof of our particular case for the reader's
convenience.

Note that for a contraction $T$ on a Hilbert space and a vector
$\xi$, the equality $\|T\xi\|=\|\xi\|$ holds if and only if
$T^*T\xi=\xi$. This shows that (ii) is equivalent to (iv).
Clearly, then (iii) also is equivalent to (iv).

For any $x\in N_1$, we have $\theta(x^*)\theta(x)\le\theta(x^*x)$
by Schwarz inequality.
Since $\nu_2$ is faithful, it follows that the equality
$||\theta(x)||_2=||x||_2$ holds if and only if
$\theta(x^*)\theta(x)=\theta(x^*x)$. It is well-known that (i) is
equivalent to the conditions $\theta(x^*)\theta(x)=\theta(x^*x)$
and $\theta(x)\theta(x^*)=\theta(xx^*)$. Thus (i) implies (ii) and
(iii), and (ii) and (iii) together imply (i). Since we already
know that (ii) and (iii) are equivalent, we conclude that all four
conditions are equivalent. \ep

For $s\in I=\Irr(G)$, consider the map $\Theta_s\colon
L^\infty(G)\to B(H_s)$,
$$
\Theta_s(a)=(\varphi\otimes\iota)(U^s(a\otimes 1){U^s}^*)
=\Theta(a)I_s.
$$
The following lemma shows that there exists a map
$\Theta^*_s\colon B(H_s)\to L^\infty(G)$ such that
$$
\phi_s(\Theta_s(a)x)=\varphi(a\Theta^*_s(x)),
\ \ a\in L^\infty(G), x\in B(H_s).
$$

\begin{lem}
We have
$$
\Theta^*_s(x)=(\iota\otimes\omega_s)\Phi_l(x)
=(\iota\otimes\omega_s)({U^s}^*(1\otimes x)U^s).
$$
\end{lem}

\bp Recall that for the modular group $\sigma^\varphi_t$ we have
\begin{equation} \label{e2.10}
(\sigma^\varphi_t\otimes\iota)(W)
=(1\otimes\rho^{it})W(1\otimes\rho^{it}).
\end{equation}
On the other hand, $\sigma^{\phi_s}_t(x)=\rho^{-it}x\rho^{it}$ for
$x\in B(H_s)$. Thus
$(\sigma^\varphi_t\otimes\sigma^{\phi_s}_t)(U^s)
=U^s(1\otimes\rho^{2it})$. Hence
\begin{eqnarray*}
\phi_s(\Theta_s(a)x)&=&(\varphi\otimes\phi_s)(U^s(a\otimes1)
{U^s}^*(1\otimes
x))\\&=&(\varphi\otimes\phi_s)((a\otimes1){U^s}^*(1\otimes
x)(\sigma^\varphi_{-i}\otimes\sigma^{\phi_s}_{-i})(U^s))\\
&=&(\varphi\otimes\phi_s(\cdot\rho^2))((a\otimes1){U^s}^*(1\otimes
x)U^s),
\end{eqnarray*}
whence $\Theta^*_s(x)=(\iota\otimes\omega_s)({U^s}^*(1\otimes
x)U^s)$, because $\omega_s=\phi_s(\cdot\rho^2)$. \ep

As a byproduct we get $\phi_s\Theta_s=\varphi$. Note also that
both maps $\Theta_s$ and $\Theta^*_s$ are $G$-equivariant.

Now we can compute $\Theta^*$.

\begin{lem}
Let $\phi\in\C_l(\hat G)$ be a generating state. Set
$\nu_0=\hat\eps|_{H^\infty(G,\phi)}$. Then for the Poisson
integral $\Theta\colon (L^\infty(G),\varphi)\to(H^\infty(\hat
G,\phi),\nu_0)$, we have
$$
\Theta^*(x)=s^*-\lim_{n\to\infty}\sum_{s\in
I}\phi^n(I_s)\Theta^*_s(x).
$$
\end{lem}

\bp First note that under the identification of $H^\infty(\hat
G,\phi)$ with the relative commutant for a product-type action,
the state $\nu_0$ coincides with the product-state. It follows
that $\nu_0$ is faithful, and its modular group is given by the
restriction of the modular group
$\sigma^{\hat\psi}_t=\Ad\rho^{-it}$ of the right-invariant Haar
weight $\hat\psi$ on $l^\infty(\hat G)$ to $H^\infty(\hat
G,\phi)$. It is then easy to see that $\nu_0\Theta=\varphi$ and
using (\ref{e2.10}) that
$\sigma^{\nu_0}_t\Theta=\Theta\sigma^\varphi_t$. Hence $\Theta^*$
indeed exists.

Recall \cite[Theorem 3.6(2)]{I} that as $\phi$ is generating, the
product on $H^\infty(\hat G,\phi)$ is given by
$$
x\cdot y=s^*-\lim_{n\to\infty}P^n_\phi(xy).
$$
Thus for any $a\in L^\infty(G)$ and $x\in H^\infty(\hat G,\phi)$
we have
\begin{eqnarray*}
\nu_0(\Theta(a)\cdot x)&=&\lim_{n\to\infty}\hat\eps
P^n_\phi(\Theta(a)x)
=\lim_{n\to\infty}\phi^n(\Theta(a)x)\\
&=&\lim_{n\to\infty}\sum_{s\in I}\phi^n(I_s)\phi_s(\Theta_s(a)x)\\
&=&\lim_{n\to\infty}\sum_{s\in
I}\phi^n(I_s)\varphi(a\Theta^*_s(x)),
\end{eqnarray*}
so  $\sum_{s\in I}\phi^n(I_s)\Theta^*_s(x)\to\Theta^*(x)$ in
weak operator topology. Note, however, that by
$G$-equivariance of $\Theta^*_s$ if $x$ is in the spectral
subspace of $H^\infty(\hat G,\phi)$ corresponding to an
irreducible representation of $G$, then $\Theta^*_s(x)$ is in the
spectral subspace of $L^\infty(G)$, which is finite dimensional.
It follows that the convergence is in norm on a dense
$*$-subalgebra of $H^\infty(\hat G,\phi)$. Since
$$
\varphi\left(\sum_s\phi^n(I_s)\Theta^*_s(x)\right)
=\sum_s\phi^n(I_s)\phi_s(x)=\phi^n(x)=\nu_0(x)
$$
for $x\in H^\infty(\hat G,\phi)$, the operators
$\sum_s\phi^n(I_s)\Theta^*_s$ are contractions with respect to the
$L^2$-norms. Hence we have $s^*$-convergence on the whole space
$H^\infty(\hat G,\phi)$.\ep

{\it From now onwards we assume that the counit $\eps$ is bounded
on $C(G)$. This is the case for $q$-deformations of compact
connected semisimple Lie groups.}

\medskip

Since the map $\Theta^*\Theta$ is $G$-equivariant, it maps the
spectral subspaces of $L^\infty(G)$ into themselves. The same is
true for $\Theta^*_s\Theta_s$. It follows that for any $a\in C(G)$
the sequence $\{\sum_s\phi^n(I_s)\Theta^*_s\Theta_s(a)\}_n$ is in
$C(G)$, and it converges in norm to $\Theta^*\Theta(a)$. Note now
that $G$-equivariance of $\Theta^*\Theta$ implies that
$$
\Theta^*\Theta(a)=(\iota\otimes\eps)\D\Theta^*\Theta(a)
=(\iota\otimes\eps\Theta^*\Theta)\D(a).
$$
So to prove that $\Theta^*\Theta(a)=a$ for an element $a\in C(G)$,
it is enough to show that $\eps\Theta^*\Theta(b)=\eps(b)$ for any
element $b$ of the form $(\omega\otimes\iota)\D(a)$, $\omega\in
C(G)^*$. As $\eps\Theta^*_s=\omega_s$ and
$$
\sum_s\phi^n(I_s)\omega_s=\sum_s\phi^n(I_s)\phi_s(\cdot\rho^2)
=\phi^n(\cdot\rho^2)=\omega^n,
$$
where $\omega=\phi(\cdot\rho^2)\in\C_r(\hat G)$, we have
$$
\eps\Theta^*\Theta(b)=\lim_{n\to\infty}\omega^n\Theta(b).
$$
Thus using equivalence of (i) and (iv) in Lemma \ref{2.1} we get
the following criterion for multiplicativity of $\Theta$.

\begin{prop} \label{2.4}
Let $\phi\in\C_l(\hat G)$ be a generating state. Set
$\omega=\phi(\cdot\rho^2)$. Then the sequence
$\{\omega^n\Theta\}^\infty_{n=1}$ of states on $C(G)$ is
$w^*$-convergent. For a $G$-invariant subspace $X$ of $C(G)$ (that
is, $\D(X)\subset C(G)\otimes X$), the limit state coincides with
the counit $\eps$ on $X$ if and only if $X$ is in the
multiplicative domain of the Poisson integral $\Theta\colon
L^\infty(G)\to H^\infty(\hat G,\phi)$.\epp
\end{prop}

The operators $\Theta_s$ and $\Theta^*_s$ are analogues of
well-known classical constructions~\cite{Be,Per}. Let for the
moment $G$ be a compact Lie group, $U\colon G\to B(H)$ a finite
dimensional unitary representation. Fix a vector $\xi\in H$,
$\|\xi\|=1$. Let $T\subset G$ be the stabilizer of the line
$\7C\xi$. For an operator $S\in B(H)$, its covariant Berezin
symbol $\sigma(S)$ is defined by $\sigma(S)(g)=(SU_g\xi,U_g\xi)$.
The covariant symbol $\sigma$ is a $G$-equivariant map from $B(H)$
into $C(G/T)$. Consider the inner products on $C(G/T)$ and $B(H)$
given by the $G$-invariant probability measure and the normalized
trace, respectively. Then there exists an adjoint
$\breve\sigma\colon C(G/T)\to B(H)$ of $\sigma$. Explicitly,
$$
\breve\sigma(f)=d\int f(g)U_gP_\xi U^*_gdg,
$$
where $d=\dim H$ and $P_\xi$ is the projection onto $\7C\xi$. A
function $f$ is called a contravariant Berezin symbol of
$\breve\sigma(f)$. The map $B=\sigma\breve\sigma$ is called the
Berezin transform.

If we consider $U$ as a corepresentation of $C(G)$, then the
definition of $\sigma$ can be written as
$$
\sigma(S)=(\iota\otimes\omega_\xi)(U^*(1\otimes S)U),
$$
where $\omega_\xi=(\cdot\xi,\xi)$ is the vector-state defined by
$\xi$. Thus we see that our operator $\Theta^*_s$ is just $\sigma$
with $\omega_\xi$ replaced by $\omega_s$. Then
$\Theta_s=(\Theta^*_s)^*$ is an analogue of $\breve\sigma$.

Suppose now that $G$ is a semisimple Lie group, $U=U^\lambda\colon
G\to B(H_\lambda)$ an irreducible representation with highest
weight $\lambda$, $\xi=\xi_\lambda$ a highest weight vector. Let
$B_\lambda$ be the corresponding Berezin transform.
Note that $\omega^n_{\xi_\lambda}=\omega_{\xi_{n\lambda}}$.
It is proved
in~\cite{D} that the sequence $\{B_{n\lambda}\}^\infty_{n=1}$
converges to the identity on $C(G/T)$ as $n\to\infty$. This is a
key step to show that the full matrix algebras $B(H_{n\lambda})$,
$n\in\7N$, provide a quantization of $C(G/T)$, see~\cite{L,R}. In
view of the $G$-equivariance of the Berezin transform, it is
enough to prove the convergence at the unit of~$G$. The proof is
based on the observation that the states $\eps B_{n\lambda}$ are
given by measures which are absolutely continuous with respect to
the Haar measure and such that the Radon-Nikodym derivatives,
up to normalization, are powers of a single function $h$ such that
$h(g)=1$ for $g\in T$ and $h(g)<1$ for $g\notin T$. The proof of
our $q$-analogue will be based on the study of ergodic properties
of an auxiliary operator.

\medskip

For a normal linear functional $\omega$ on $l^\infty(\hat G)$
define a linear operator $A_\omega\colon C(G)\to C(G)$ by
$$
A_\omega(a)=(\iota\otimes\omega)\hat\Phi(a)
=(\iota\otimes\omega)(W(a\otimes1)W^*).
$$
Then $\omega\Theta=\varphi A_\omega$. Since $\hat\Phi$ is a right
action of $\hat G$, so that
$(\iota\otimes\hat\D)\hat\Phi=(\hat\Phi\otimes\iota)\hat\Phi$, we
have $A_{\omega_1\omega_2}=A_{\omega_1}A_{\omega_2}$ for any
$\omega_1,\omega_2\in l^\infty(\hat G)_*$. Thus
$\omega^n\Theta=\varphi A^n_\omega$, and by Proposition~\ref{2.4}
to show multiplicativity of $\Theta$ on the quantum flag manifold,
it is enough to prove the following result.

\begin{thm} \label{2.5}
Let $G=SU_q(n)$ ($0<q<1$), $T\subset SU_q(n)$ the maximal torus,
and $\omega\in\C_r(\hat G)$ a normal right $G$-invariant state,
$\omega\ne\hat\eps$. Then the counit $\eps$ is the only
$A_\omega$-invariant state on $C(G/T)$.
\end{thm}

We will prove the result by induction. For this we first have to
establish functorial properties of the operators $A_\omega$.

Let $H$ be a closed subgroup of $G$. By this we mean that $H$ a
compact quantum group and that we are given a surjective unital
$*$-homomorphism $\pi\colon C(G)\to C(H)$ which respects
comultiplication. We can define a left and a right action of $H$
on $C(G)$ by the homomorphisms
$$
C(G)\to C(H)\otimes C(G), \ \ a\to(\pi\otimes\iota)\D,
$$
and
$$
C(G)\to C(G)\otimes C(H), \ \ a\to(\iota\otimes\pi)\D,
$$
respectively. The corresponding fixed point algebras are denoted
by $C(H\backslash G)$ and $C(G/H)$.

By considering the elements of $l^\infty(\hat G)$ as linear
functionals on $\A(G)$ we can define a dual homomorphism
$\hat\pi\colon l^\infty(\hat H)\to l^\infty(\hat G)$.
Equivalently, one can consider the unitary corepresentation
$U=(\pi\otimes\iota)(W)$ of $C(H)$, where $W$ is the multiplicative
unitary for $G$, and set $\hat\pi=\pi_U$.

\begin{lem} \label{2.6}
Let $H$ be a closed subgroup of $G$ defined by
$\pi\colon C(G)\to C(H)$. Then \enu{i} $\pi
A_\omega=A_{\omega\hat\pi}\pi$ for any $\omega\in l^\infty(\hat
G)_*$; \enu{ii} $A_\omega(C(G/H))\subset C(G/H)$ for any
$\omega\in l^\infty(\hat G)_*$;\enu{iii} $A_\omega(C(H\backslash
G))\subset C(H\backslash G)$ for any $\omega\in\C_r(\hat G)$.
\end{lem}

\bp Let $W$ and $W_0$ be the multiplicative unitaries for $G$ and
$H$, respectively. Set $U=(\pi\otimes\iota)(W)$ as above, so that
$\hat\pi=\pi_U$. As
$U=(\iota\otimes\pi_U)(W_0)=(\iota\otimes\hat\pi)(W_0)$, we get
$$
\pi A_\omega(a)=(\iota\otimes\omega)(U(\pi(a)\otimes1)U^*)
=(\iota\otimes\omega\hat\pi)(W_0(\pi(a)\otimes1)W_0^*)
=A_{\omega\hat\pi}\pi(a),
$$
which proves (i).

\smallskip

Let $a\in C(G/H)$, so that $(\iota\otimes\pi)\D(a)=a\otimes1$.
Then using the pentagon equation $W_{12}W_{13}W_{23}=W_{23}W_{12}$
we get
\begin{eqnarray*}
(\iota\otimes\pi)\D A_\omega(a) &=&(\iota\otimes\pi\otimes\omega)
(W^*_{12}W_{23}(1\otimes a\otimes1)W^*_{23}W_{12})\\
&=&(\iota\otimes\pi\otimes\omega)(W_{13}W_{23}W^*_{12}(1\otimes
a\otimes1)W_{12}W^*_{23}W^*_{13})\\
&=&(\iota\otimes\pi\otimes\omega)
(W_{13}W_{23}(\D(a)\otimes1)W^*_{23}W^*_{13})\\
&=&(\iota\otimes\pi\otimes\omega)
(W_{13}W_{23}(a\otimes1\otimes1)W^*_{23}W^*_{13})\\
&=&A_\omega(a)\otimes1,
\end{eqnarray*}
which shows (ii).

\smallskip

Suppose now that $a\in C(H\backslash G)$. Then
\begin{eqnarray*}
(\pi\otimes\iota)\D A_{\omega}(a)
&=&(\pi\otimes\iota\otimes\omega)
(W^*_{12}W_{23}(1\otimes a\otimes1)W^*_{23}W_{12})\\
&=&(\pi\otimes\iota\otimes\omega) (W_{13}W_{23}W^*_{12}(1\otimes
a\otimes1)W_{12}W^*_{23}W^*_{13})\\
&=&(\pi\otimes\iota\otimes\omega)
(W_{13}W_{23}(\D(a)\otimes1)W^*_{23}W^*_{13})\\
&=&(\pi\otimes\iota\otimes\omega) (W_{13}W_{23}(1\otimes
a\otimes1)W^*_{23}W^*_{13})\\
&=&(\pi\otimes\iota\otimes\omega) (W_{23}(1\otimes
a\otimes1)W^*_{23})\\
&=&1\otimes A_{\omega}(a),
\end{eqnarray*}
where in the next to last equality we used right $G$-invariance of
$\omega$. This proves (iii).\ep

\begin{lem} \label{2.7}
For any $\omega\in\C_r(\hat G)$, we have
$\omega\hat\pi\in\C_r(\hat H)$.
\end{lem}

An equivalent way of saying this is that if $\rho_0=f_1$ is the
Woronowicz character for $H$, then $\hat\pi(\rho_0)\rho^{-1}$
commutes with $\hat\pi(l^\infty(\hat H))$. Yet another equivalent
statement is that the homomorphism~$\pi$ intertwines the scaling
groups of $C(G)$ and $C(H)$.

\medskip

\bpp{Lemma \ref{2.7}} Keeping the notation of the proof of the
previous lemma, consider the right actions $\alpha=\alpha_{W,r}$
and $\alpha_0=\alpha_{U,r}$ of $G$ and $H$, respectively, on
$B(L^2(G))$. Extend $\omega$ to a normal $G$-invariant functional
$\tilde\omega$ on $B(L^2(G))$. Since $(\pi\otimes\iota)(W)=U$, we
have $(\iota\otimes\pi)\alpha(x)=\alpha_0(x)$ for any $x\in
B(L^2(G))$ (note that the expression $(\iota\otimes\pi)\alpha(x)$
makes sense as we have a well-defined homomorphism
$\iota\otimes\pi\colon M(B_0(L^2(G))\otimes C(G))\to
M(B_0(L^2(G))\otimes C(H))$). Hence $\tilde\omega$ is
$H$-invariant, so $\omega\hat\pi=\tilde\omega\pi_U\in\C_r(\hat
H)$. \ep

We can now lay the foundation for our induction argument.

\begin{lem}
Suppose $\eta$ is a state on $C(G)$ such that $\eta=\eps$ on
$C(H\backslash G/H)=C(H\backslash G)\cap C(G/H)$. Then there
exists a state $\eta_0$ on $C(H)$ such that $\eta_0\pi=\eta$.
\end{lem}

\bp Suppose $a\ge0$, $\pi(a)=0$. We have to prove that
$\eta(a)=0$. Let $\varphi_0$ be the Haar state on $C(H)$. We have
$$
(\varphi_0\pi\otimes\iota\otimes\varphi_0\pi)\D^2(a)\in
C(H\backslash G/H),
$$
where $\D^2=(\D\otimes\iota)\D=(\iota\otimes\D)\D$. Hence
\begin{eqnarray*}
(\varphi_0\pi\otimes\eta\otimes\varphi_0\pi)\D^2(a)
&=&(\varphi_0\pi\otimes\eps\otimes\varphi_0\pi)\D^2(a)
=(\varphi_0\pi\otimes\varphi_0\pi)\D(a)\\
&=&(\varphi_0\otimes\varphi_0)\D_0\pi(a)=0,
\end{eqnarray*}
where $\D_0$ is the comultiplication on $C(H)$. Since the state
$\varphi_0\otimes\varphi_0$ is faithful by our assumptions on
quantum groups, we conclude that
$$
(\pi\otimes\eta\otimes\pi)\D^2(a)=0.
$$
Applying $\eps_0\otimes\iota\otimes\eps_0$, where $\eps_0$ is the
counit on $C(H)$, and using $\eps_0\pi=\eps$ we get
$\eta(a)=0$.\ep

\begin{cor} \label{2.9}
Let $T\subset H\subset G$ be compact quantum groups, and
$\pi\colon C(G)\to C(H)$ the homomorphism defining the inclusion
$H\hookrightarrow G$. Let $\omega$ be a state in $\C_r(\hat G)$.
Assume that \enu{i} the counit $\eps$ on $C(G)$ is the only
$A_\omega$-invariant state on $C(H\backslash G/H)$; \enu{ii} the
counit $\eps_0$ on $C(H)$ is the only
$A_{\omega\hat\pi}$-invariant state on $C(H/T)$.

Then $\eps$ is the only $A_\omega$-invariant state on $C(G/T)$.
\end{cor}

\bp If $\eta$ is $A_{\omega}$-invariant, then $\eta=\eps$ on
$C(H\backslash G/H)$, so by the previous lemma $\eta=\eta_0\pi$
for some state $\eta_0$ on $C(H)$. By Lemma~\ref{2.6}(i) and
surjectivity of $\pi$, the state $\eta_0$ is
$A_{\omega\hat\pi}$-invariant. Hence $\eta_0=\eps_0$ on $C(H/T)$.
As $\pi(C(G/T))\subset C(H/T)$, we have
$\eta=\eta_0\pi=\eps_0\pi=\eps$ on $C(G/T)$. \ep

Let us now recall the definition of $SU_q(n)$,
see e.g.~\cite{KS}.

The algebra $C(U_q(n))$ of continuous functions on the compact
quantum group $U_q(n)$ is generated by $n^2+1$ elements $u_{ij}$,
$1\le i,j\le n$, $t$ satisfying the relations
$$
u_{ik}u_{jk}=qu_{jk}u_{ik},\ \
u_{ki}u_{kj}=qu_{kj}u_{ki}\ \ \hbox{for}\ i<j,
$$
$$
u_{il}u_{jk}=u_{jk}u_{il}\ \ \hbox{for}\ i<j,\ k<l,
$$
$$
u_{ik}u_{jl}-u_{jl}u_{ik}=(q-q^{-1})u_{jk}u_{il} \ \
\hbox{for}\ i<j,\ k<l,
$$
$$
\det_q(U)t=t\det_q(U)=1,\ \ u_{ij}t=tu_{ij}\ \ \hbox{for any}\
i,j,
$$
where $U=(u_{ij})_{i,j}$ and $\det_q(U)=\sum_{w\in
S_n}(-q)^{\ell(w)}u_{w(1)1}\ldots u_{w(n)n}$, with $\ell(w)$ being
the number of inversions in $w\in S_n$. The involution on
$C(U_q(n))$ is given by
$$
t^*=\det_q(U),\ \
u^*_{ij}=(-q)^{j-i}\det_q(U^{\hat i}_{\hat j})t,
$$
where $U^{\hat i}_{\hat j}$ is the matrix obtained from $U$ by
removing the $i$th row and the $j$th column. Taking the quotient
of $C(U_q(n))$ by the ideal generated by $\det_q(U)-1$, we obtain
the algebra $C(SU_q(n))$.

If $m<n$, then $U_q(m)\times\7T^{n-m}$ is a subgroup of $U_q(n)$.
Namely, if $u'_{ij}$, $1\le i,j\le m$, and $t'$ are the
generators of $C(U_q(m))$, and $z_1,\ldots,z_{n-m}$ are the
canonical generators of $C(\7T^{n-m})$, then the homomorphism
$C(U_q(n))\to C(U_q(m)\times\7T^{n-m})$ is given by
$$
u_{ij}\mapsto\cases{u'_{ij} & if $1\le i,j\le m$, \cr z_{i-m} & if
$i=j>m$, \cr 0 & otherwise,}
$$
and thus $t\mapsto t'z^{-1}_1\ldots z^{-1}_{n-m}$. Taking the
intersection with $SU_q(n)$, in other words, taking the quotient
of $C(U_q(m)\times\7T^{n-m})$ by the ideal generated by
$1-t'z^{-1}_1\ldots z^{-1}_{n-m}$, we get a subgroup
$S(U_q(m)\times\7T^{n-m})\cong U_q(m)\times\7T^{n-m-1}$ of
$SU_q(n)$. The subgroup $T=S(\7T^n)\cong\7T^{n-1}$ we call the
maximal torus of~$SU_q(n)$.

Consider now the filtration $T=G_1\subset G_2\subset\ldots\subset
G_n=SU_q(n)$ of $SU_q(n)$, where $G_m=S(U_q(m)\times\7T^{n-m})$,
and let $\pi_m\colon C(SU_q(n))\to C(G_m)$ be the corresponding
homomorphisms. We will prove by induction that the counit on
$C(G_m)$ is the only $A_{\omega\hat\pi_m}$-invariant state on
$C(G_m/T)$. For $m=1$ there is nothing to prove as $G_1=T$. Thus
by Corollary~\ref{2.9} we just have to show that the counit is the
only $A_{\omega\hat\pi_{m+1}}$-invariant state on $C(G_m\backslash
G_{m+1}/G_m)$.

Consider the embedding $SU_q(2)\hookrightarrow SU_q(n)$
corresponding to the right lower corner of $SU_q(m+1)$. In other
words, if $u'_{ij}$, $1\le i,j\le 2$, are the generators of
$C(SU_q(2))$, we define a homomorphism $\theta_m\colon
C(SU_q(n))\to C(SU_q(2))$ by
$$
\theta_m(u_{ij})\mapsto\cases{u'_{i-m+1,j-m+1} & if $m\le i,j\le
m+1$, \cr \delta_{ij} & otherwise.}
$$
Note that $\theta_m$ factorizes through $\pi_{m+1}\colon
C(SU_q(n))\to C(G_{m+1})$, so $\theta_m=\theta'_m\pi_{m+1}$
for some homomorphism $\theta'_m\colon C(G_{m+1})\to C(SU_q(2))$.

\begin{lem} \label{2.10}
For each $m$, $1\le m\le n-1$, the homomorphism $\theta'_m$ maps
$C(G_m\backslash G_{m+1}/G_m)$ isomorphically onto
$C(\7T\backslash SU_q(2)/\7T)$.
\end{lem}

\bp Note that $G_m\backslash G_{m+1}/G_m\cong
(U_q(m)\times\7T)\backslash U_q(m+1)/(U_q(m)\times\7T)$. Then the
result can be deduced from Theorem 4.7 in \cite{MNY}, which
implies that $C((U_q(m)\times\7T)\backslash
U_q(m+1)/(U_q(m)\times\7T))$ is generated by
$u_{m+1,m+1}u^*_{m+1,m+1}$.

Another possibility is to use the classification of irreducible
representations of algebras of functions on homogeneous spaces,
see e.g.~\cite{PV,DS}.  Namely, let either $G=U_q(n)$,
$H=U_q(n-1)\times\7T$ and $T=\7T^n$, or $G=SU_q(n)$,
$H=S(U_q(n-1)\times\7T)$ and $T=S(\7T^n)\cong\7T^{n-1}$. Then the
irreducible representations of $C(G)$ can be described as
follows~\cite{KoS}. Consider the homomorphisms $\theta_k\colon
C(G)\to C(SU_q(2))$, $1\le k\le n-1$, and $\pi_1\colon C(G)\to
C(T)$ as above. Fix an irreducible
infinite dimensional representation $\pi$ of $C(SU_q(2))$.
Then for a character
$\chi$ of $T$ and an element $w\in S_n$ with a reduced
decomposition $w=\tau_{i_1}\ldots\tau_{i_k}$, where
$\tau_i=(i,i+1)$, $1\le i\le n-1$, we set
$$
\pi_{w,\chi}=(\pi\theta_{i_1})\times\ldots\times(\pi\theta_{i_k})
\times(\chi\pi_1)
=((\pi\theta_{i_1})\otimes\ldots\otimes(\pi\theta_{i_k})
\otimes(\chi\pi_1))\D^k.
$$
Up to equivalence the representation $\pi_{w,\chi}$ is independent
of the reduced decomposition of $w$, and $\{\pi_{w,\chi}\}_{w\in
S_n,\chi\in\hat T}$ is a complete set of irreducible
representations of $C(G)$. If $w\in S_{n-1}$, then the
representation $\pi_{w,\chi}$ factorizes through $C(H)$, so its
restrictions to $C(G/H)$ and $C(H\backslash G)$ are given by the
counit. If $w\in S_n\backslash S_{n-1}$, then $w$ can be written
as $w_1\tau_{n-1}w_2$ with $\ell(w)=\ell(w_1)+\ell(w_2)+1$ and
$w_1,w_2\in S_{n-1}$. Since the representations $\pi_{w_1,0}$ and
$\pi_{w_2,\chi}$ factorize through $C(H)$, and
$$
\D^2(C(H\backslash G/H))\subset C(H\backslash G)\otimes C(G)
\otimes C(G/H),
$$
we see that $\pi_{w,\chi}(a)=1\otimes(\pi\theta_{n-1})(a)\otimes1$
for any $a\in C(H\backslash G/H)$. Thus the restriction of any
irreducible representation of $C(G)$ to $C(H\backslash G/H)$
factorizes through $\theta_{n-1}\colon C(G)\to C(SU_q(2))$. Hence
$\theta_{n-1}\colon C(H\backslash G/H)\to C(SU_q(2))$ is
injective. As $T\subset H$, the image is obviously contained in
$C(\7T\backslash SU_q(2)/\7T)$. In fact, it coincides with
$C(\7T\backslash SU_q(2)/\7T)$, since e.g. it is easy to check
that $\theta_{n-1}(u_{nn}u^*_{nn})$ generates $C(\7T\backslash
SU_q(2)/\7T)$. \ep

Now to prove Theorem \ref{2.5} we just have to show that the
counit on $C(SU_q(2))$ is the only
$A_{\omega\hat\theta_m}$-invariant state on $C(\7T\backslash
SU_q(2)/\7T)$ for $1\le m\le n-1$. Note that as $SU(n)$ is a
simple Lie group, the restriction of a non-trivial representation
of $SU(n)$ to a non-discrete subgroup is non-trivial. It follows
that if $\omega$ is non-trivial on $l^\infty(\widehat{SU_q(n)})$,
that is, $\omega\ne\hat\eps$, then $\omega\hat\theta_m$ is
non-trivial on $l^\infty(\widehat{SU_q(2)})$. Thus we can reduce
the proof of Theorem~\ref{2.5} to the case of $SU_q(2)$, moreover,
in this case it suffices to prove that the counit is the only
$A_{\omega}$-invariant state on $C(\7T\backslash SU_q(2)/\7T)$.
For this we could use the results of~\cite{I,NT} saying that the
Poisson integral is a homomorphism on $C(SU_q(2)/\7T)$, at least
under some summability assumptions on $\omega$. We will instead
give a self-contained probabilistic proof.

\smallskip

Let $\{u_{ij}\}_{1\le i,j\le2}$ be the generators of $C(SU_q(2))$.
Set $\alpha=u_{11}$ and $\gamma=u_{21}$. Then
$$
U=\pmatrix{u_{11} & u_{12}\cr u_{21} & u_{22}}
=\pmatrix{\alpha & -q\gamma^*\cr \gamma & \alpha^*},
$$
and the relations can be written as
$$
\alpha^*\alpha+\gamma^*\gamma =1,
\ \ \alpha\alpha^* +q^2\gamma^*\gamma =1,\ \
\gamma^*\gamma =\gamma\gamma^*,\ \
\alpha\gamma =q\gamma\alpha,\ \ \alpha\gamma^*=q\gamma^*\alpha.
$$
The comultiplication $\D$ is determined by the formulas
$$
\D (\alpha )=u_{11}\otimes u_{11}+u_{12}\otimes u_{21}
=\alpha\otimes\alpha -q\gamma^*\otimes\gamma,\ \
\D (\gamma )=\gamma\otimes\alpha +\alpha^*\otimes\gamma.
$$
The homomorphism $C(SU_q(2))\to C(\7T)$ is given by $\alpha\mapsto
z$, $\gamma\mapsto0$. The monomials $\alpha^k(\gamma^*)^l\gamma^m$
and $(\alpha^*)^k(\gamma^*)^l\gamma^m$, $k,l,m\ge0$, span a dense
$*$-subalgebra of $C(SU_q(2))$. It is then easy to see that
$C(\7T\backslash SU_q(2)/\7T)$ is generated by $\gamma^*\gamma$.
The spectrum of $\gamma^*\gamma$ is the set
$I_{q^2}=\{0\}\cup\{q^{2n}\}^\infty_{n=0}$. Thus we can identify
$C(\7T\backslash SU_q(2)/\7T)$ with the algebra $C(I_{q^2})$ of
continuous functions on $I_{q^2}$. Under this identification the
counit is given by the evaluation at $0\in I_{q^2}$. The Markov
operator $A_\omega$ defines a random walk on
$I_{q^2}\backslash\{0\}$. If this random walk is transient, then
$\nu A^n_\omega\to\eps$ as $n\to\infty$ for any state $\nu$ on
$C(\7T\backslash SU_q(2)/\7T)$. In particular, the counit $\eps$
is the only $A_\omega$-invariant state. As was remarked in
\cite{NT}, transience of a random walk on a non-Kac discrete
quantum group follows easily from the fact that the Markov
operator has a positive eigenvector with eigenvalue strictly
smaller than $1$. It is natural to expect that the operator
$A_\omega$ acting on the dual side has the same property.

\begin{prop} \label{2.11}
Consider the function $f$ on $I_{q^2}$ defined by $f(0)=0$,
$f(q^{2k})=a_k$, where $\{a_k\}^\infty_{k=0}$ is the sequence
defined by the recurrence relation
$$
a_0=1,
$$
$$
2(1-q^{2k+1})a_k=q^{-1}(1-q^{2k+2})a_{k+1}+q(1-q^{2k})a_{k-1}.
$$
Then \enu{i} $f\in C(I_{q^2})$ and $f(q^{2k})>0$ for any $k\ge0$;
\enu{ii} for any normal right $SU_q(2)$-invariant state
$\omega=\sum_s\lambda_s\omega_s\in\C_r(\widehat{SU_q(2)})$, the
element $f$ is an eigenvector for $A_\omega$ with eigenvalue
$\displaystyle\sum_s\lambda_s\frac{2s+1}{[2s+1]_q}$.
\end{prop}

As usual, we identify the set $\Irr(SU_q(2))$ with $\2\7Z_+$. Then
$\displaystyle [2s+1]_q=\frac{q^{2s+1}-q^{-2s-1}}{q-q^{-1}}$ is
the quantum dimension of the representation with spin
$s\in\2\7Z_+$.

\medskip

\bpp{Proposition \ref{2.11}} The proof of (i) is analogous to that
of \cite[Lemma~5.4]{I}. To see that $a_k>0$, rewrite the
recurrence relation as
$$
q^{-1}(1-q^{2k+2})(a_{k+1}-qa_k)
=(1-q^{2k})(a_k-qa_{k-1})+q^{2k}(1-q)^2a_k.
$$
It follows by induction that $a_{k+1}-qa_k\ge0$, so $a_k\ge q^k$.
It remains to show that $a_k\to0$. It is clear that the sequence
$\{a_k\}^\infty_{k=0}$ cannot grow faster than a geometric
progression. Hence the generating function
$$
g(z)=\sum^\infty_{k=0}a_kz^k
$$
is analytic in a neighborhood of zero. The recurrence relation can
then be written as
$$
2(g(z)-qg(q^2z))=q^{-1}z^{-1}(g(z)-g(q^2z))+qz(g(z)-q^2g(q^2z)),
$$
that is,
$$
g(z)=\left(\frac{1-q^2z}{1-qz}\right)^2g(q^2z).
$$
We see that $g$ extends to a meromorphic function with poles at
$z=q^{-2k-1}$, $k\ge0$. In particular, the series $\sum_ka_kz^k$
converges for $|z|<q^{-1}$, whence $a_k\to0$. In fact, since
$$
\lim_{z\to q^{-1}}(1-qz)^2g(z)
=\prod^\infty_{k=0}\left(\frac{1-q^{2k+1}}{1-q^{2k+2}}\right)^2
=\frac{(q;q^2)^2_\infty}{(q^2;q^2)^2_\infty},
$$
we have $\displaystyle
a_k\sim kq^k\frac{(q;q^2)^2_\infty}{(q^2;q^2)^2_\infty}$.

\medskip

To prove (ii), first consider the case $\omega=\omega_{\2}$, so that
$$
\omega=\frac{1}{[2]_q}\Tr\left(\cdot\pmatrix{q^{-1} & 0\cr 0 &
q}\right)
$$
as $U^\2=U$. Then we have
\begin{eqnarray*}
A_\omega(x)&=& \frac{1}{[2]_q}\Tr\left(\pmatrix{\alpha &
-q\gamma^*\cr \gamma & \alpha^*} \pmatrix{x & 0\cr 0 &
x}\pmatrix{\alpha^* & \gamma^*\cr -q\gamma &
\alpha}\pmatrix{q^{-1} & 0\cr 0 &
q}\right)\\
&=& \frac{1}{[2]_q}(q^{-1}(\alpha
x\alpha^*+q^2\gamma^*x\gamma)+q(\gamma
x\gamma^*+\alpha^*x\alpha)).
\end{eqnarray*}
Identifying $C(\7T\backslash SU_q(2)/\7T)=C^*(\gamma^*\gamma)$
with $C(I_{q^2})$, and using the identities
$$
\alpha^*(\gamma^*\gamma)^k\alpha=q^{-2k}(\gamma^*\gamma)^k
(1-\gamma^*\gamma),\ \
\alpha(\gamma^*\gamma)^k\alpha^*=q^{2k}(\gamma^*\gamma)^k
(1-q^2\gamma^*\gamma),
$$
we see that the action of $A_\omega$ on the functions on $I_{q^2}$
is given by
$$
(A_\omega h)(t)=\frac{1}{[2]_q}
\Big(q^{-1}\Big((1-q^2t)h(q^2t)+q^2th(t)\Big)
 +q\Big(th(t)+(1-t)h(q^{-2}t\Big)\Big)\Big).
$$
Then the definition of $\{a_k\}_k$ shows that
$A_{\omega}f=\frac{2}{[2]_q}f$. To prove that $f$ is an eigenvalue
for $A_{\omega_s}$ for any $s$, recall from \cite[Section~6]{I}
that the identity
$$
\omega_s\omega_\2=\frac{d_{s-\2}}{d_sd_\2}\omega_{s-\2}
+\frac{d_{s+\2}}{d_sd_\2}\omega_{s+\2},
$$
where $d_s=[2s+1]_q$, implies that there exists a polynomial
$p_{2s}$ of degree $2s$ such that $p_{2s}(\omega_\2)=\omega_s$.
Then $A_{\omega_s}=p_{2s}(A_{\omega_\2})$, so $f$ is an
eigenvector for $A_{\omega_s}$ with eigenvalue
$p_{2s}(\frac{2}{[2]_q})$. As
$\frac{2s+1}{[2s+1]_q}=\omega_s(\rho^{-1})$, and $\rho$ is
group-like, we have
$$
p_{2s}\left(\frac{2}{[2]_q}\right)=p_{2s}(\omega_\2(\rho^{-1}))
=p_{2s}(\omega_\2)(\rho^{-1})=\omega_s(\rho^{-1})
 =\frac{2s+1}{[2s+1]_q},
$$
which finishes the proof of Proposition.
\ep

It follows that the random walk defined by $A_\omega$ on
$I_{q^2}\backslash\{0\}$ is transient for any normal
right-invariant state $\omega\ne\hat\eps$. More precisely, if
$A_\omega f=\lambda f$, then the probability of visiting a point
$t_2$ from a point $t_1$ at the $n$th step is not larger than
$f(t_1)f(t_2)^{-1}\lambda^n$. Hence $\nu
A^n_\omega\to\delta_0=\eps$ as $n\to\infty$ for any state $\nu$ on
$C(I_{q^2})$. Note that to see that $\eps$ is the only
$A_\omega$-invariant state is even easier. Indeed, if $\nu$ is
$A_\omega$-invariant, we have $\nu(f)=\lambda\nu(f)$, so
$\nu(f)=0$ and $\nu=\delta_0$. This completes the proof of
Theorem~\ref{2.5}.

\bigskip

\section{Random walk on the center} \label{s3}

To prove surjectivity of the Poisson integral, we will obtain an
estimate on the dimensions of the spectral subspaces of
$H^\infty(\hat G,\phi)$. By a result of Hayashi~\cite{H}, if the
fusion algebra of a group $G$ is commutative, then any central
harmonic element is a scalar. Equivalently, the action of $G$ on
the Poisson boundary is ergodic. This already implies that the
spectral subspaces of $H^\infty(\hat G,\phi)$ are finite
dimensional, more precisely, the dimension of the spectral
subspace corresponding to an irreducible representation $U$ is not
larger than the square of the quantum dimension of
$U$~\cite{B,HLS}. This estimate is clearly not sufficient for our
purposes. We will show that in our situation ergodicity of the
action provides a better estimate.

\smallskip

The result of Hayashi is in fact more general. It was obtained as
a consequence of an analogue of double ergodicity of the Poisson
boundary, see e.g. \cite{K}. Since in our situation the proof can
be made more concrete, we will present a detailed argument.

For Markov operators $P$ and $Q$ on a von Neumann algebra, we say
that an element $x$ is $(P,Q)$-harmonic if $P(x)=Q(x)=x$. Let now
$\phi$ and $\omega$ be generating normal states on $l^\infty(\hat
G)$ with the same support. Set
$$
(N,\nu)=\Big(\mathop{\otimes}^{-1}_{-\infty}(l^\infty(\hat
G),\phi)\Big)\bigotimes\Big(\mathop{\otimes}^{+\infty}_0
(l^\infty(\hat G),\omega)\Big),
$$
and let $\gamma$ be the shift to the right on $N$. For any finite
interval $I=[n,m]\subset\7Z$ we have a normal homomorphism
$j_I\colon l^\infty(\hat G)\to N$ defined by $\hat\D^{(m-n)}$.
Then the space of $(P_\phi,Q_\omega)$-harmonic elements, where
$P_\phi=(\phi\otimes\iota)\hat\D$ and
$Q_\omega=(\iota\otimes\omega)\hat\D$, can be embedded in the
space $N^\gamma$ of $\gamma$-invariant elements by the
homomorphism $j_\7Z$,
$$
j_\7Z(x)=s^*-\lim_{n\to-\infty\atop m\to+\infty}j_{[n,m]}(x).
$$
Note that if $\phi\ne\omega$, then the automorphism $\gamma$ is
never ergodic~\cite[Lemma~4]{K}. Nevertheless
the following result holds.

\begin{prop}{\cite[Proposition~3.4]{H}}
Let $\phi\in\C_l(\hat G)$ be a generating normal left
$G$-invariant state on $l^\infty(\hat G)$,
$\omega=\phi(\cdot\rho^2)\in\C_r(\hat G)$ the corresponding
right-invariant state. Then the space of central
$(P_\phi,Q_\omega)$-harmonic elements consists of the scalars.
\end{prop}

\bp Note that by definition $\phi^n=\omega^n$ on the center of
$l^\infty(\hat G)$. Consider the operator
$Q_\phi=(\iota\otimes\phi)\hat\D$. Although in general $Q_\phi\ne
Q_\omega$, for any left-invariant state $\phi'$ we have
$$
\phi'Q^n_\phi=\phi^n P_{\phi'}=\omega^n P_{\phi'}=\phi'Q^n_\omega
$$
on the center. Similarly, for any right-invariant state $\omega'$
we have $\omega'P^n_\omega=\omega'P^n_\phi$ on the center. It
follows that for any central elements $z_1$ and $z_2$ and any
intervals $I_1\subset I_2$, $I_1=[n,m]$, $I_2=[n-k,m+l]$, we have
$$
\nu(j_{I_1}(z_1)j_{I_2}(z_2))=\phi^{m-n+1}(z_1 Q^l_\omega
P^k_\phi(z_2)).
$$
Indeed, e.g. in the case when $m+l<0$ we get
$$
\nu(j_{I_1}(z_1)j_{I_2}(z_2))=\phi^{m-n+1}(z_1 Q^l_\phi
P^k_\phi(z_2))=\phi^{m-n+1}(z_1 Q^l_\omega P^k_\phi(z_2)).
$$
Hence for any central $z_1$ and $z_2$ and any finite intervals
$I_1\subset I_2$ we have
$$
||j_{I_1}(z_1)-j_{I_2}(z_2)||_2 =||\gamma j_{I_1}(z_1)-\gamma
j_{I_2}(z_2)||_2
$$
Thus if $z$ is a central $(P_\phi,Q_\omega)$-harmonic element,
then the distance $||j_\7Z(z)-\gamma^nj_I(z)||_2$ is independent
of $n$. Since on the one hand this distance goes to zero as
$I\nearrow\7Z$, and on the other hand $\gamma^nj_I(z)$ converges
in weak operator topology to a scalar as $n\to+\infty$, we
conclude that $j_\7Z(z)$ is a scalar. \ep

\begin{cor}{\cite[Corollary~3.5]{H}}
Assume that the fusion algebra $R(G)$ of the group $G$ is
commutative. Then for any generating state $\phi\in\C_l(\hat G)$,
the scalars are the only central $P_\phi$-harmonic elements.
\end{cor}

\bp Commutativity of the fusion algebra means that
$P_\phi=Q_\omega$ on the center. Thus any central
$P_\phi$-harmonic element is $(P_\phi,Q_\omega)$-harmonic, and we
can apply the previous proposition. \ep

\medskip

Consider now the random walk on the center in more detail.
Identify the center of $l^\infty(\hat G)$ with $l^\infty(I)$,
where $I=\Irr(G)$. For a fixed generating state $\phi\in\C_l(\hat
G)$, let $\{p(s,t)\}_{s,t\in I}$ be the transition probabilities
defined by the restriction of $P_\phi$ to $l^\infty(I)$, so
$P_\phi(I_t)I_s=p(s,t)I_s$. Let $(\Omega,\7P_0)$ be the path space
of the corresponding random walk,
$$
\Omega=\prod^\infty_{n=1}I,\ \
\7P_0(\{\underline{s}\,|\,s_1=t_1,\ldots,s_n=t_n\})
=p(0,t_1)p(t_1,t_2)\ldots p(t_{n-1},t_n).
$$
Denote by $\pi_n$ the projection $\Omega\to I$ onto the $n$th
factor.

Similarly to Section~\ref{s1}, set
$$
(N,\nu)=\mathop{\otimes}^{-1}_{-\infty}(l^\infty(\hat G),\phi),
$$
$j_n(x)=\ldots\otimes1\otimes\hat\D^{n-1}(x)$ for $x\in
l^\infty(\hat G)$, and $j_\infty(x)=s^*-\lim_n j_n(x)$ for $x\in
H^\infty(\hat G,\phi)$. (In Section~\ref{s1} we
embedded $F_n(l^\infty(\hat G))$ into $B(H)$ for some $H$ and
extended $\phi$ to a $G$-invariant normal faithful state
$\tilde\phi$ on $B(H)$, which we don't do now as the relative
commutant interpretation of the Poisson boundary will not be
important.) As was remarked in~\cite{NT}, there is an embedding
$$
j^\infty\colon(L^\infty(\Omega,\7P_0),\7P_0)\hookrightarrow(N,\nu)
$$
such that $f\pi_n\mapsto j_n(f)$ for any $f\in l^\infty(I)\subset
l^\infty(\hat G)$. If $f\in l^\infty(I)$ is harmonic, then the
sequence $\{f\pi_n\}^\infty_{n=1}$ is a martingale, so it
converges almost everywhere (a.e.) to a function $f_\infty\in
L^\infty(\Omega,\7P_0)$. Then $j_\infty(f)=j^\infty(f_\infty)$.

Denote by $H^\infty(I,\phi)$ the space $l^\infty(I)\cap
H^\infty(\hat G,\phi)$ of central harmonic elements. Let $E\colon
l^\infty(\hat G)\to l^\infty(I)$ be the unique $G$-equivariant
conditional expectation,
$$
E(x)=(\varphi\otimes\iota)\Phi_l(x)=\sum_{s\in I}\phi_s(x)I_s.
$$
By restricting $E$ to $H^\infty(\hat G,\phi)$ we get a conditional
expectation $H^\infty(\hat G,\phi)\to H^\infty(I,\phi)$.

\begin{prop} \label{3.3}
Let $x,y\in H^\infty(\hat G,\phi)$. Then the sequence
$\{f_n\}^\infty_{n=1}$ of functions on $\Omega$ defined by
$$
f_n(\underline{s})=\phi_{s_n}(xy)
$$
converges a.e. to $E(x\cdot y)_\infty\in L^\infty(\Omega,\7P_0)$.
\end{prop}

Since $f_n=E(xy)\pi_n$, one can equivalently state that
$\{E(xy)\pi_n\}_n$ and $\{E(x\cdot y)\pi_n\}_n$ converge a.e. to
the same limit.

\medskip

\bpp{Proposition \ref{3.3}}
Let $\alpha$ be the product-type action of $G$ on $N$, and
$\tilde E=(\varphi\otimes\iota)\alpha\colon N\to N^\alpha$ the
$\nu$-preserving conditional expectation. Since $j_n(x)\to
j_\infty(x)$ and $j_n(y)\to j_\infty(y)$ in $s^*$-topology, we
have
$$
j_n(xy)=j_n(x)j_n(y)\to j_\infty(x)j_\infty(y)=j_\infty(x\cdot y)
$$
in $s^*$-topology. Hence $\tilde Ej_n(xy)\to\tilde
Ej_\infty(x\cdot y)$, and as $\tilde Ej_n=j_nE$, we get
$j_nE(xy)\to j_\infty E(x\cdot y)$ in strong$^*$ operator
topology. Using that $f_n=E(xy)\pi_n$,
$j_nE(xy)=j^\infty(E(xy)\pi_n)$ and $j_\infty E(x\cdot
y)=j^\infty(E(x\cdot y)_\infty)$, we conclude that $f_n\to
E(x\cdot y)_\infty$ in measure. It remains to show that the
sequence $\{f_n\}_n$ is a.e. convergent.

It is enough to consider the case $x=y^*$. Let
$L^\infty(I^n,\7P^{(n)}_0)$ be the subalgebra of
$L^\infty(\Omega,\7P_0)$ consisting of the functions depending
only on the first $n$ coordinates, $E_n\colon
L^\infty(\Omega,\7P_0)\to L^\infty(I^n,\7P^{(n)}_0)$ the
$\7P_0$-preserving conditional expectation. For any $f\in
l^\infty(I)$ we have $E_n(f\pi_{n+1})=P_\phi(f)\pi_n$. As
$$
y^*y=P_\phi(y)^*P_\phi(y)\le P_\phi(y^*y)
$$
by Schwarz inequality, we have $E(y^*y)\le EP_\phi(y^*y)=P_\phi
E(y^*y)$, whence
$$
f_n=E(y^*y)\pi_n\le P_\phi(E(y^*y))\pi_n=E_n(E(y^*y)\pi_{n+1})
=E_n(f_{n+1}).
$$
Thus the sequence $\{f_n\}^\infty_{n=1}$ is a bounded
submartingale. By Doob's theorem, see e.g.~\cite{KSK}, it must
converge a.e. \ep

\begin{cor} \label{3.4}
Let $\phi\in\C_l(\hat G)$ be a generating state. Assume that the
Poisson boundary of the center is trivial, i.e.
$H^\infty(I,\phi)=\7C1$. Then for any $x,y\in H^\infty(\hat
G,\phi)$ and almost every path $\underline{s}\in\Omega$, we have
$\phi_{s_n}(xy)\to\hat\eps(x\cdot y)$ as $n\to\infty$.\epp
\end{cor}

\begin{cor}
Let $\phi\in\C_l(\hat G)$ be a generating state, $V$ an
irreducible representation of $G$. Assume that the Poisson
boundary of the center is trivial. Then the multiplicity of $V$ in
$H^\infty(\hat G,\phi)$ is not larger than the supremum of the
multiplicities of $V$ in $\overline{U}\times U$ for all
irreducible representations $U$ of $G$.
\end{cor}

\bp By the previous corollary, for any finite dimensional subspace
$X$ of $H^\infty(\hat G,\phi)$ and almost every path
$\underline{s}\in\Omega$, the restrictions of the irreducible
representations $l^\infty(\hat G)\to B(H_{s_n})$ to $X$ are
asymptotically isometric in $L^2$-norm as $n\to\infty$. In
particular, these restrictions are eventually injective. Since the
maps $H^\infty(\hat G,\phi)\to B(H_s)$ are $G$-equivariant, it
follows that the multiplicity of $V$ in $H^\infty(\hat G,\phi)$ is
not larger than the supremum of the multiplicities of $V$ in
$B(H_U)$ for all irreducible representations $U$ of $G$ on $H_U$.
It remains to note that the $G$-module $B(H_U)$, or more
precisely, the $\A(\hat G)$-module such that $\omega x=(\hat
S(\omega)\otimes\iota)\alpha_{U,l}(x)$ for $\omega\in\A(\hat G)$
and $x\in B(H_U)$, is isomorphic to $\overline{H}_U\otimes
H_U$.\ep

For the $q$-deformation $G$ of a compact connected semisimple Lie
group the last estimate is optimal. Indeed, let $T\subset G$ be
the maximal torus. Then for an irreducible representation $V$ of
$G$ on $H$, the multiplicity of $V$ in $L^\infty(G/T)$ is equal to
the dimension of the space of zero weight vectors in $H$, that is,
the space of $T$-invariant vectors. On the other hand, by the
Frobenius reciprocity the multiplicity of $V$ in
$\overline{U}\times U$ is the same as the multiplicity of $U$ in
$U\times V$. Both the multiplicity $N^U_{U,V}$ of $U$ in $U\times
V$ and the dimension $m_0(V)$ of the space of zero weight vectors
are known to be independent of the deformation parameter. Hence
$N^U_{U,V}\le m_0(V)$, see e.g.~\cite[\S 131]{Z}. It follows that
the spectral subspaces of $H^\infty(\hat G,\phi)$ are not larger
than the spectral subspaces of $L^\infty(G/T)$. Note also that as
$\hat\eps\Theta=\varphi$ and $\varphi$ is faithful, the Poisson
integral is injective on its multiplicative domain. Thus we get
the following result.

\begin{thm} \label{3.6}
Let $G$ be the $q$-deformation of a compact connected semisimple
Lie group, $T\subset G$ the maximal torus, $\phi\in\C_l(\hat G)$ a
generating state. Assume that the Poisson integral $\Theta\colon
L^\infty(G/T)\to H^\infty(\hat G,\phi)$ is a homomorphism. Then it
is an isomorphism.\epp
\end{thm}

Theorems A and B now follow from Proposition~\ref{2.4},
Theorem~\ref{2.5} and Theorem~\ref{3.6}. Indeed, it follows
immediately that if $\phi\in\C_l(l^\infty(\widehat{SU_q(n)}))$ is
a generating state, then the Poisson integral $\Theta\colon
L^\infty(SU_q(n)/\7T^{n-1})\to H^\infty(\widehat{SU_q(n)},\phi)$
is a $SU_q(n)$- and $\widehat{SU_q(n)}$-equivariant isomorphism,
where $\7T^{n-1}=S(\7T^n)$ is the maximal torus in $SU_q(n)$. If
$\phi$ is not generating, then $\cup_k\supp\phi^k$ corresponds to
a quotient $SU(n)/\Gamma$ of $SU(n)$ and to the quotient
$G=SU_q(n)/\Gamma$ of $SU_q(n)$, which we call the $q$-deformation
of $SU(n)/\Gamma$. More explicitly, if $\Gamma$ is the group of
roots of unity of order $m$, $m|n$, then $C(G)$ is the subalgebra
of $C(SU_q(n))$ generated by the matrix coefficients of $U^{\times
m}$, where $U$ is the fundamental representation of $SU_q(n)$. Set
$T=\7T^{n-1}/\Gamma$, so $C(T)$ is the image of $C(G)$ under the
homomorphism $C(SU_q(n))\to C(\7T^{n-1})$. Then
$L^\infty(G/T)=L^\infty(SU_q(n)/\7T^{n-1})\subset
L^\infty(SU_q(n))$. Since the assumptions of Theorem~\ref{2.5}
don't require $\omega=\phi(\cdot\rho^2)$ to be generating on
$\widehat{SU_q(n)}$, we again conclude that the Poisson integral
$\Theta\colon L^\infty(G/T)\to H^\infty(\hat G,\phi)$ is a $G$-
and $\hat G$-equivariant isomorphism. This completes the proof of
Theorem~B. To prove Theorem~A, note that the fixed point algebra
is independent of whether we consider the action of $SU_q(n)$, or
the action of its quotient $G$ such that
$\Irr(G)=\cup_k\supp\phi^k$. Thus $(N^\alpha)'\cap N$ is
$G$-equivariantly isomorphic to
$L^\infty(G/T)=L^\infty(SU_q(n)/\7T^{n-1})$. Clearly, the
isomorphism is $SU_q(n)$-equivariant. If we identify
$(N^\alpha)'\cap N$ with $H^\infty(\hat G,\phi)$, so we get an
action of $\hat G$ on $(N^\alpha)'\cap N$, then the isomorphism is
also $\hat G$-equivariant.

\bigskip

\section{Concluding remarks}

For any non-trivial product-type action of $SU_q(n)$ on $N$ the
fixed point algebra $N^{SU_q(n)}$ is obviously strictly contained
in the fixed point algebra $N^T$ for the action of the maximal
torus $T\subset SU_q(n)$ (contrary to what is claimed
in~\cite{S}). Moreover, by our results
$$
(N^{SU_q(n)})'\cap N^T\cong L^\infty(T\backslash SU_q(n)/T).
$$
Consider the product-type action defined by the fundamental
representation of $SU_q(n)$ on $\7C^n$. Let $\H_\infty(q)$ be the
Hecke algebra, that is, the algebra with generators
$g_1,g_2,\ldots$ and relations
$$
g_i^2=(q-q^{-1})g_i+1,\ \
g_ig_{i+1}g_i=g_{i+1}g_ig_{i+1},\ \
g_ig_j=g_jg_i\ \ \hbox{for}\ \ |i-j|\ge2.
$$
Then $N^{SU_q(n)}$ is the weak operator closure of the image of
$\H_\infty(q)$ under the homomorphism $\pi\colon\H_\infty(q)\to N$
defined by
$$
\pi(g_1)=\ldots\otimes1\otimes\left(q\sum_im_{ii}\otimes m_{ii}
+(q-q^{-1})\sum_{i<j}m_{ii}\otimes m_{jj}+\sum_{i\ne j}m_{ij}
\otimes m_{ji}\right),
$$
$\pi(g_n)=\gamma^{n-1}\pi(g_1)$, where $\gamma\colon N\to N$ is
the shift to the left, see e.g. \cite[Proposition 8.40]{KS}. The
homomorphism $\pi$ should not be confused with the
homomorphisms $\pi_+$ and $\pi_-$ defined by
$$
\pi_\pm(g_1)=\ldots\otimes1\otimes\left(q\sum_im_{ii}\otimes m_{ii}
+(q-q^{-1})\sum_{i>j}m_{ii}\otimes m_{jj}\pm\sum_{i\ne j}m_{ij}
\otimes m_{ji}\right),
$$
$\pi_\pm(g_n)=\gamma^{n-1}\pi_\pm(g_1)$. Using that the unique
left $SU_q(n)$-invariant state on $B(\7C^n)$ is defined by the
density matrix
$$
\frac{1-q^2}{1-q^{2n}}\pmatrix{q^{2(n-1)} & 0 &\ldots & 0\cr
0 & q^{2(n-2)} & \ldots & 0\cr
\vdots & \vdots & \ddots & \vdots \cr
0 & 0 & \ldots & 1},
$$
we see that if $E\colon N\to\gamma(N)$ is the $\nu$-preserving
conditional expectation, then $E\pi_+(g_1)$ and $E\pi_-(g_1)$ are
scalars, while $E\pi(g_1)$ is not. According to~\cite{PP,S1}, one
has $\pi_+(\H_\infty(q))''=\pi_-(\H_\infty(q))''=N^T$.

\medskip

As we showed in Section~\ref{s2}, if an element $a\in C(G)$ is in
the multiplicative domain of $\Theta$, then
$$
a=\lim_{n\to\infty}\sum_{s\in I}\phi^n(I_s)\Theta^*_s\Theta_s(a).
$$
If in addition the Poisson boundary of the center is trivial, then
we have a stronger convergence result:
\begin{equation} \label{e4.1}
a=\lim_{n\to\infty}\Theta^*_{s_n}\Theta_{s_n}(a)
\end{equation}
for almost every path $\underline{s}\in\Omega$. Indeed, first note
that as we assume that $l^\infty(\hat G)$ has a generating state,
the space $C(G)$ is separable. Hence by Corollary~\ref{3.4} we
have $\phi_{s_n}(\Theta(a)\Theta(b))\to
\hat\eps(\Theta(a)\cdot\Theta(b))$ for almost every path
$\underline{s}\in\Omega$ and any $b\in C(G)$. Since
$$
\phi_{s_n}(\Theta(a)\Theta(b))
=\varphi(\Theta^*_{s_n}\Theta_{s_n}(a)b)\ \ \hbox{and}\ \
\hat\eps(\Theta(a)\cdot\Theta(b))=\hat\eps(\Theta(ab))=\varphi(ab),
$$
we see that convergence~(\ref{e4.1}) holds in weak operator
topology for almost every path $\underline{s}\in\Omega$. As in
Section~\ref{s2}, by $G$-equivariance of
$\Theta^*_{s_n}\Theta_{s_n}$ we conclude that the convergence is
in norm.

Let now $G=SU_q(n)$. Identify $I=\Irr(SU_q(n))$ with the set of
dominant weights of $SU(n)$. Let $V_\lambda$ be an irreducible
representation of $SU_q(n)$ with highest weight $\lambda$. What we
used in Section~\ref{s3}, is that the multiplicity of $V_\lambda$
in $B(H_s)$ is not larger than the multiplicity of $V_\lambda$ in
$C(SU_q(n)/T)$. In fact, the multiplicities are equal as soon as
$s$ is sufficiently large, see~\cite{Z}. It is e.g. enough to
require $s+w\lambda$ to be dominant for any element $w$ of the
Weyl group. Recall also that in the classical case Berezin
transforms converge to the identity on the flag manifold along any
ray in the Weyl chamber. Thus it is natural to conjecture that
convergence~(\ref{e4.1}) holds for every sequence
$\underline{s}=\{s_n\}^\infty_{n=1}$ such that the distance from
$s_n$ to the walls of the Weyl chamber goes to infinity.

\medskip

A significant part of our results is valid for $q$-deformations of
arbitrary compact connected semisimple Lie groups. The point where
we crucially used that the group was $SU_q(n)$, was
Lemma~\ref{2.10}, which allowed us to reduce the proof of
Theorem~\ref{2.5} to the study of a one-dimensional random walk.
Lemma~\ref{2.10} is also valid for $q=1$, in which case it is an
immediate consequence of the fact that $S(U(m)\times\7T)\subset
SU(m+1)$ is a Riemannian symmetric pair of rank one. Hence there
is hope that similar considerations could work for $SO(n)$,
$Sp(n)$ and $F_4$, see~\cite[Ch.~X, Table~V]{He}. This will be
discussed in detail elsewhere. For the exceptional groups $E_6$,
$E_7$, $E_8$ and $G_2$, however, our reduction procedure leads us
to consider random walks of higher dimensions. The ultimate goal
would of course be to find a unified proof. For this it could be
instructive to understand the origin of the eigenvector
constructed in Proposition~\ref{2.11}, since to prove that $\eps$
is the only $A_\omega$-invariant state on $C(G/T)$, it is enough
to find a strictly positive eigenvector for $A_\omega$ in the
kernel of $\eps$ on $C(G/T)$ with eigenvalue less than $1$. Remark
also that by using commutativity of the fusion algebra, one can
show that if $G$ is the $q$-deformation of a compact connected
simple Lie group and $\eps$ is the only $A_\omega$-invariant state
on $C(G/T)$ for some $\omega\ne\hat\eps$, then $\eps$ is the only
invariant state for any $\omega\ne\hat\eps$.

\bigskip

\flushleft{Masaki Izumi, Department of Mathematics, Graduate
School of Science, Kyoto University, Sakyo-ku, Kyoto 606-8502,
Japan\\
{\it e-mail}: izumi@kusm.kyoto-u.ac.jp}

\flushleft{Sergey Neshveyev, Mathematics Institute, University of
Oslo, PB 1053 Blindern, Oslo 0316, Norway\\
{\it e-mail}: sergeyn@math.uio.no}

\flushleft{Lars Tuset, Mathematics Institute, University of
Oslo, PB 1053 Blindern, Oslo 0316, Norway\\
{\it e-mail}: larstu@math.uio.no}

\end{document}